\documentclass[a4paper,11pt]{article}
\usepackage{ifstyle}   
\usepackage{graphicx}
\usepackage{epsfig} 
\usepackage{color}
\numberwithin{equation}{section}  

\usepackage[francais]{babel}               
\usepackage[T1]{fontenc}
\usepackage[latin1]{inputenc}  
\usepackage{amssymb,diagrams}             
\usepackage{array}                 
\usepackage{fullpage}

\newtheorem{lem}[subsubsection]{Lemme}

\newcommand{\R}{\mathbb R}
\newcommand{\C}{\mathbb C}

\newcommand{\I}{\mathbb I}

\newcommand{\contract}{\mathrel{\kern-1.5pt\vrule width6.0pt height0.4pt depth0pt
                \vrule width0.4pt height4.0pt depth0pt}}
\newcommand{\retract}{\mathrel{\kern-1.5pt \vrule width0.4pt height5.0pt depth0pt
                \vrule width6.0pt height0.4pt depth0pt }}

\def\opn#1#2{\def#1{\operatorname{#2} } } 
\opn\Ric{Ricci} 
\opn\Trac{Trace} 
\opn\det{det} 
\opn\Ker{Ker} 
\opn\exp{exp}
\opn\exph{exph}
\opn\Herm{Herm} 
 
\begin{document}
{\def\thefootnote{\relax}
\footnote{\hskip-0.6cm
{\bf{Mots-cl\'es}} : Connexions de Chern, Courbure de Chern, Variétés presque complexes, Coordonnées presque complexes \\
{\bf{Classification AMS}} : 32C35.}} 
\begin{center} 
\Huge{\bf{La connexion et la courbure de Chern du fibré tangent d'une variété presque complexe}}
\\
\vspace{0.4cm}
\huge{Nefton Pali}
\end{center} 
\vspace{\fill}
{\bf Résumé}.-Sur une variété presque complexe $(X,J)$ l'opérateur $\bar{\partial}_{_{J}}$ induit une connexion de type $(0,1)$ sur le fibré des $(p,0)$-formes. Dans le cas d'une structure presque complexe intégrable cette connexion induit la structure holomorphe canonique du fibré des $(p,0)$-formes. En considérant le cas $p=1$ on peut étendre la connexion correspondante à toutes les puissances de Schur du fibré des $(1,0)$-formes. En utilisant l'isomorphisme $\C$-linéaire entre le fibré des $(1,0)$-formes et le fibré cotangent complexe $T^*_{X,J}$ on déduit aussi des connexions canoniques de type $(0,1)$ sur les puissances de Schur du fibré cotangent complexe $T^*_{X,J}$. 
\\
Dans le cas complexe intégrable ces connexions donnent les structures holomorphes canoniques de ces fibrés. Dans le cas presque complexe non intégrable les connexions en question
donnent seulement les structures holomorphes canoniques sur les restrictions des fibrés correspondants aux images des courbes $J$-holomorphes lisses. 
\\
Nous introduisons la notion de courbure de Chern pour ces fibrés, notion dont le sens géométrique est la généralisation naturelle de la notion classique de courbure de Chern pour les fibrés holomorphes sur une variété complexe.  
\\
Nous portons un intérêt particulier au cas du fibré tangent en vue des applications concernant la régularisation des fonctions $J$-plurisousharmoniques à l'aide du flot géodésique d'une connexion de Chern sur le fibré tangent (voir \cite{Pal}). Cette méthode à été déjà utilisée par Demailly \cite{Dem-2} dans le cas complexe intégrable. 
\\
Nous montrons une formule explicite qui relie la connexion de Chern du fibré tangent avec la connexion de Levi-Civita à l'aide des obstructions géométriques dérivant de la torsion de la structure presque complexe et du défaut de la métrique à être symplectique. En particulier nous donnons une formule explicite qui permet de relier la torsion de la connexion de Chern du fibré tangent avec les obstructions précédentes. Une formule qui relie les deux connexions précédentes peut être aussi trouvée dans l'article de Gauduchon \cite{Gau}. L'utilité de la connexion de Chern dans le problème de régularisation des fonctions $J$-plurisousharmoniques dérive du fait que son expression locale par rapport à des repères du fibré des $(1,0)$-vecteurs tangents est la plus simple possible parmi les connexions hermitiennes. 
\\
Ensuite nous introduisons la notion de coordonnées presque complexes au voisinage d'un point. Cette notion nous permet d'étudier la façon dont la torsion de la structure presque complexe et le caractère non symplectique de la métrique se traduisent en une obstruction à l'existence de coordonnées géodésiques complexes, qui n'existent que dans le cas Kählerien. Cette étude est nécessaire pour le calcul asymptotique du flot géodésique induit par une connexion de Chern sur le fibré tangent.
\\
\\
{\bf{Abstract}}.-The $\bar{\partial}_{_{J}}$ operator over an almost complex manifold induces canonical connections of type $(0,1)$ over the bundles of $(p,0)$-forms. If the almost complex structure is integrable then the previous connections induce the canonical holomorphic structures of the bundles of $(p,0)$-forms. For $p=1$ we can extend the corresponding connection to all Schur powers of the bundle of $(1,0)$-forms. Moreover using the canonical $\C$-linear isomorphism betwen the bundle of $(1,0)$-forms and the complex cotangent bundle $T^*_{X,J}$ we deduce canonical connections of type $(0,1)$ over the Schur powers of the complex cotangent bundle $T^*_{X,J}$. If the almost complex structure is integrable then the previous $(0,1)$-connections induces the canonical holomorphic structures of those bundles. In the non integrable case those $(0,1)$-connections induces just the holomorphic canonical structures of the restrictions of the corresponding bundles to the images of smooth $J$-holomorphic curves. We introduce the notion of Chern curvature for those bundles. The geometrical meaning of this notion is a natural generalisation of the classical notion of Chern curvature for the holomophic vector bundles over a complex manifold. We have a particular interest for the case of the tangent bundle in view of applications concerning the regularisation of $J$-plurisubharmonic fonctions by means of the geodesic flow induced by a Chern connection on the tangent bundle. This method has been used by Demailly \cite{Dem-2} in the complex integrable case. Our specific study in the case of the tangent bundle gives an asymptotic expanson of the Chern flow which relates in a optimal way the geometric obstructions caused by the torsion of the almost complex structure, and the non symplectic nature of the metric.
\section{Connexions sur les faisceaux de modules de fonctions $\ci$ au dessus des variétés presque complexes}\label{Con-presq-compl} 
Soit $(X,J)$ une variété presque complexe de classe ${\cal C}^{\infty}$ et de dimension réelle $2n$. On désigne par ${\cal E}_X\equiv{\cal E}_X(\R)$ le faisceau des fonctions ${\cal C}^{\infty}$ à valeurs réelles, par $\pi^{1,0} _{_J}:T_X\otimes_{_{\R}}\C\longrightarrow T^{1,0}_{X,J}$ la projection sur le fibré des $(1,0)$-vecteurs tangents et par $\pi^{0,1} _{_J}$ celle sur le fibré des $(0,1)$-vecteurs tangents. On désigne par $T_{X,J}$ le fibré tangent dont les fibres sont munies de la structure complexe donnée par $J$ et par
$$
{\cal E}^{p,q}_{_{X,J} }\equiv {\cal E}(\Lambda ^{p,q}_{_J}T_X^*),
\;\Lambda ^{p,q}_{_J}T_X^*:=\Lambda^p _{_{\C}}(T^{1,0}_{X,J}  )^*\otimes_{_{\C}} \Lambda^q _{_{\C}}(T^{0,1}_{X,J})^*
$$ 
le faisceau des $(p,q)$-formes par rapport à la structure presque complexe $J$. On rappelle que sur une variété presque complexe la différentielle se décompose sous la forme 
$$
d=\partial_{_J }+\bar{\partial}_{_J }-\theta_{_J }-\bar{\theta }_{_J },
$$ 
où pour toute $k$-forme complexe 
$\omega \in {\cal E}(\Lambda^k _{_{\C}}(T_X\otimes_{_{\R}}\C)^*)(U)$ au dessus d'un ouvert $U$ et tout champ de vecteurs complexes $\xi _0,...,\xi _k\in{\cal E}(T_X\otimes_{_{\R}}\C)(U)$ on a les expressions suivantes:
\begin{eqnarray*} 
&\displaystyle{\partial_{_J }\omega\, (\xi _0,...,\xi _k):=\sum_{0\leq j \leq k}(-1)^j\xi^{1,0}_j .\,\omega (\xi _0,...,\widehat{\xi _j},..., \xi _k)+ }&
\\
\\
&\displaystyle{+\sum_{0\leq j<l \leq k}(-1)^{j+l}\omega ([\xi^{1,0}_j,\xi^{1,0}_l]^{1,0}+[\xi^{0,1} _j,\xi^{1,0} _l]^{0,1}+[\xi^{1,0} _j,\xi^{0,1} _l]^{0,1},\xi _0,...,\widehat{\xi _j},...,\widehat{\xi _l},..., \xi _k)   }&
\\
\\
\\
&\displaystyle{\bar{\partial}_{_J }\omega\, (\xi _0,...,\xi _k):=\sum_{0\leq j \leq k}(-1)^j\xi^{0,1}_j .\,\omega (\xi _0,...,\widehat{\xi _j},..., \xi _k)+ }&
\\
\\
&\displaystyle{+\sum_{0\leq j<l \leq k}(-1)^{j+l}\omega ([\xi^{0,1}_j,\xi^{0,1}_l]^{0,1}+[\xi^{0,1} _j,\xi^{1,0} _l]^{1,0}+[\xi^{1,0} _j,\xi^{0,1} _l]^{1,0},\xi _0,...,\widehat{\xi _j},...,\widehat{\xi _l},..., \xi _k)   }&
\end{eqnarray*}
\begin{eqnarray*} 
&\displaystyle{\theta_{_J }\omega\, (\xi _0,...,\xi _k):=-\sum_{0\leq j<l \leq k}(-1)^{j+l}\omega ([\xi^{1,0}_j,\xi^{1,0}_l]^{0,1},\xi _0,...,\widehat{\xi _j},...,\widehat{\xi _l},..., \xi _k)   }&
\\
\\
\\
&\displaystyle{\bar{\theta }_{_J }\omega\, (\xi _0,...,\xi _k):=-\sum_{0\leq j<l \leq k}(-1)^{j+l}\omega ([\xi^{0,1}_j,\xi^{0,1}_l]^{1,0},\xi _0,...,\widehat{\xi _j},...,\widehat{\xi _l},..., \xi _k)   }&
\end{eqnarray*}
avec $\xi ^{1,0}:=\pi ^{1,0} _{_J}(\xi ),\;[\cdot,\cdot]^{1,0}:=\pi^{1,0}_{_J}[\cdot,\cdot] $ et de façon analogue pour les indices $(0,1)$. Les bidegrés des opérateurs 
$\partial_{_J },\,\bar{\partial}_{_J },\,\theta_{_J }$ et $\bar{\theta }_{_J }$ sont respectivement $(1,0),\;(0,1),\;(2,-1)$ et $(-1,2)$. En effet si 
$\omega\in{\cal E}^{p,q}_{_{X,J} } (U)$ est une $(p,q)$-forme alors les $(p+q+1)$-formes $\partial_{_J }\omega,\,\bar{\partial}_{_J }\omega,\,\theta_{_J }\omega,\,\bar{\theta }_{_J }\omega$ sont nulles en restriction aux fibrés $\Lambda ^{r,s}_{_J}T_X,\,r+s=p+q+1$ respectivement aux bidegrés $(r,s)\not=(p+1,q),\,(r,s)\not=(p,q+1),\,(r,s)\not=(p+2,q-1)
,\,(r,s)\not=(p-1,p+2)$. On déduit alors que l'opérateur $T=\partial_{_J },\;\bar{\partial}_{_J },\;\theta_{_J }$ où $\bar{\theta }_{_J }$ vérifie la règle de Leibnitz 
$$
T(u\wedge v)= Tu\wedge v+(-1)^{\deg\,u} u\wedge  Tv.
$$
On a aussi les formules $\overline{(\partial_{_J }\,u)}=\bar{\partial}_{_J } \,\bar{u},\,\overline{(\theta_{_J }\,u)}=\bar{\theta }_{_J }\,\bar{u}   $.
\begin{defi}
On désigne par $\tau_{_J}\in {\cal E}(\Lambda ^{2,0}_{_J}T_X^*\otimes_{_{\C}} T^{0,1}_{X,J})(X)$ le tenseur de la torsion de la structure presque complexe définie par la formule $\tau_{_J}(\xi ,\eta):= [\xi^{1,0},\eta^{1,0}]^{0,1}$ pour tout $\xi,\eta\in{\cal E}(T_X\otimes_{_{\R}}\C)(U)$, où $U\subset X$ désigne un ouvert quelconque. Le tenseur de  la structure presque complexe est dit intégrable si $\tau_{_J}=0$.
\end{defi}
On remarque que $\tau_{_J}=0$ si et seulement si $\theta_{_J}=0$, si et seulement si $d=\partial_{_J }+\bar{\partial}_{_J }$. 
\\
\\
{\bf Note au lecteur}. Le $\C$-isomorphisme canonique $T_{X,J,x}\rightarrow T^{1,0}_{X,J,x}$ implique le $\C$-isomorphisme  
$\Lambda^{p,q}_{_J}T_{X,x} ^*\otimes_{_{\C}}T_{X,J,x}\rightarrow\Lambda^{p,q}_{_J}T_{X,x} ^*\otimes_{_{\C}} T^{1,0}_{X,J,x},\,\alpha \mapsto u$. Pour tout vecteur réel $\xi\in\Lambda^{p+q}_{_{\R} }T_{X,x} $ on a l'égalité $\alpha (\xi )=u(\xi)+\overline{u(\xi)}$. en effet soit $(\zeta _k)_k\subset (T^{1,0}_{X,J,x})^{\oplus n}$ un repère complexe de $T^{1,0}_{X,J,x}$. Alors $(v_k)_k\subset (T_{X,J,x})^{\oplus n},\,v_k=\zeta _k+\bar{\zeta }_k$ est un repère complexe de $T_{X,J,x}$. La forme $\alpha$ s'écrit alors sous la forme $\alpha =\sum_k\,\alpha _k\otimes_{_{J} }v_k,\;\alpha \in \Lambda^{p,q}_{_J}T_{X,x}^*$ et $u=\sum_k\,\alpha _k\otimes \zeta _k$. Pour tout élément $\xi \in\Lambda^{p+q}_{_{\C} }(T_X\otimes_{_{\R}} \C)$ on a par définition  
$$
\alpha(\xi ) =\sum_k\,\alpha _k(\xi )\times_{_{J} } v_k=\sum_k\,(\alpha _k(\xi )\zeta _k+\overline{\alpha (\xi)}\bar{\zeta}).
$$
Si $\xi \in\Lambda^{p+q}_{_{\R} }T_{X,x} \subset\Lambda^{p+q}_{_{\C} }(T_{X,x} \otimes_{_{\R}} \C)$ on a l'égalité voulue. Nous considérons l'espace vectoriel
$$
R^{p,q}_{_{J} } (T_{X,x} \otimes_{_{\R}} \C):=\{u+\bar{u}\;|\;u\in\Lambda ^{p,q}_{_J}T_{X,x} ^*\otimes_{_{\C}} T^{1,0} _{X,J,x} \}
$$
avec la structure de produit $\times_{_{J}}$ définie par la formule $c\,\times_{_{J} }(u+\bar{u}):=cu+\overline{cu},\,c\in \C$. 
Le fait qu'une forme $\C$-linéaire sur le complexifié $T_{X,x} \otimes_{_{\R}} \C$ de l'espace tangent $T_{X,x} $ soit déterminée de façon univoque à partir de sa restriction à $T_{X,x}$ nous suggère qu'il est très naturel de considérer le $\C$-isomorphisme
$\Lambda^{p,q}_{_J}T_{X,x} ^*\otimes_{_{\C}}T_{X,J,x}\rightarrow R^{p,q}_{_{J} } (T_{X,x} \otimes_{_{\R}} \C),\;\alpha \mapsto u+\bar{u}$. Dans la suite on identifiera donc les éléments de l'espace vectoriel 
$\Lambda ^{p,q}_{_J}T_{X,x} ^*\otimes_{_{\C}} T_{X,J,x}$ avec les éléments du type $u+\bar{u},\,u\in\Lambda ^{p,q}_{_J}T_{X,x} ^*\otimes_{_{\C}} T^{1,0} _{X,J,x}$. L'utilité d'un tel formalisme sera clarifié dans la suite.
\hfill $\Box$
\\ 
\\
On définit le tenseur de Nijenhuis 
$$
N_{_J}\in {\cal E}(\Lambda ^{0,2}_{_J}T_X^*\otimes_{_{\C}} T_{X,J})(X)
$$
par la formule $N_{_{J}}:=\tau_{_J}+\bar{\tau}_{_{J}}$. Bien évidemment $N_{_{J}}=0 $ si et seulement si $\tau_{_J}=0$. Il est élémentaire de vérifier l'identité:
$$
4N_{_{J}}(\xi ,\eta) =[\xi ,\eta]+J[\xi ,J\eta]+J[J\xi ,\eta]-[J\xi ,J\eta]
$$
pour tout champ de vecteurs complexes
$\xi ,\eta\in {\cal E}(T_X\otimes_{_{\R}}\C)(X)$. On rappelle le célèbre théorème de Newlander-Nirenberg (voir \cite{We}, \cite{Hor}, \cite{Dem-1}, chapitre VIII, \cite{mal}, \cite{Nij-Woo} et \cite{New-Nir}).
\begin{theo}{\bf (Newlander-Nirenberg)}.
Soit $(X,J)$ une variété presque complexe. L'existence d'une structure holomorphe ${\cal O}_X$ sur la variété $X$ telle que la structure presque complexe associée $J_{{\cal O}_X } $ soit égale à $J$ est équivalente à l'intégrabilité de la structure presque complexe $J$.
\end{theo} 
On considère les définitions suivantes.
\begin{defi} Soient $(X,J_1)$ et $(Y,J_2)$ deux variétés presque complexes et $f:X\rightarrow Y$ une application différentiable. L'application $f$ est dite $(J_1,J_2)$-holomorphe si sa différentielle vérifie la condition $J_2(f(x))\cdot d_xf=d_xf\cdot J_1(x)$ pour tout $x\in X$.
\end{defi}
Pour tout  application différentiable $f:X\rightarrow Y$, la différentielle 
$df\in \Gamma (X, T^*_X\otimes_{_{\R} }f^*T_Y)$ se décompose sous la forme 
$$
df=\partial_{_{J_1,J_{2 }}}f+\bar{\partial}_{_{J_{1 },J_{2 }}}f,
$$
où
\begin{eqnarray*} 
&\displaystyle{
\partial_{_{J_1,J_{2 }}}f_{_{|_x } }:=\frac{1}{2}\Big(d_xf-J_2(f(x))\cdot d_xf\cdot J_1(x)\Big) }&
\\
\\
&\displaystyle{
 \bar{\partial}_{_{J_{1 },J_{2 }}}f_{_{|_x } }:=\frac{1}{2}\Big(d_xf+J_2(f(x))\cdot d_xf\cdot J_1(x)\Big). }&
\end{eqnarray*}
Bien évidemment 
\begin{eqnarray*} 
&\displaystyle{
\partial_{_{J_1,J_{2 }}}f\in \Gamma (X, T^*_{X,J_1} \otimes_{_{\C} }f^*T_{Y,J_2} )}&
\\
\\
&\displaystyle{
\bar{\partial}_{_{J_{1 },J_{2 }}}f\in \Gamma (X, T^*_{X,-J_1} \otimes_{_{\C} }f^*T_{Y,J_2} ) }&
\end{eqnarray*}
 et l'application $f$ est $(J_1,J_2)$-holomorphe si et seulement si $\bar{\partial}_{_{J_{1 },J_{2 }}}f=0$.

\begin{defi} 
Soit $(X,J)$ une variété presque complexe et $(\Sigma,j)$ une courbe holomorphe lisse. Une courbe $(j,J)$-holomorphe est une application différentiable 
$\gamma :(\Sigma ,j) \longrightarrow (X,J)$  dont la différentielle vérifie la condition 
$J(\gamma (z))\cdot d_z\gamma=d_z\gamma \cdot j$ pour tout $z\in \Sigma$. On désigne par $i$ la structure presque complexe canonique sur $\R^2\equiv \C$. Une courbe $J$-holomorphe locale est une courbe $(i,J)$-holomorphe $\gamma :(B^1_{\delta},i) \longrightarrow (X,J)$ définie sur le disque complexe de rayon $\delta >0$.
\end{defi}  
On a alors qu'une application différentiable 
$\gamma :B^1_{\delta} \longrightarrow X$ est une courbe $J$-holomorphe locale si et seulement si elle vérifie l'équation
$\bar{\partial}_{_{j,J}}\gamma (\frac{\partial}{\partial t})=0,\;z=t+is$ qui s'écrit explicitement sous la forme
$$
\partial_s\gamma =J(\gamma )\cdot \partial_t\gamma, 
$$
où $\partial_s\gamma:=d\gamma(\frac{\partial}{\partial s})$. On peut montrer, (voir prop.2.3.6 dans l'article de Sikorav, dans l'ouvrage \cite{Ad}) que si $\gamma$ est une courbe $J$-holomorphe alors $\gamma \in \ci(B^1_{\delta };X)$.
On aura besoin aussi de la définition suivante.
\begin{defi} 
Soit ${\cal G}$ un faisceau de ${\cal E}(\C)$-modules sur $X$. Une connexion sur le faisceau ${\cal G}$ est un morphisme de faisceaux de groupes additifs $\nabla_{_{\cal G}}:{\cal G}\longrightarrow{\cal G}\otimes_{_{{\cal E}}}{\cal E}(T^*_X)\simeq{\cal G}\otimes_{_{{\cal E}(\C)}}{\cal E}(T^*_X\otimes_{_{\R}}\C) $ tel que 
$\nabla_{_{\cal G}}(g\cdot f)=\nabla_{{\cal G}}\,g\cdot f+g\otimes df$ pour tout $g\in{\cal G}(U)$ et $f\in {\cal E}(\C)(U)$, où $U\subset X$ est un ouvert quelconque.
\end{defi} 
La donnée d'une connexion $\nabla_{_{\cal G}}$ sur le faisceau de ${\cal E}(\C)$-modules ${\cal G}$ détermine de façon univoque une dérivation sur le complexe $({\cal G}\otimes_{_{{\cal E}}}{\cal E}(\Lambda ^k_{_{\R}} T^*_X))_{k\geq 0} $. En effet on peut définir l'extension 
$$
\nabla_{_{\cal G}}:{\cal G}\otimes_{_{{\cal E}}}{\cal E}(\Lambda ^k_{_{\R}}T^*_X)
\longrightarrow{\cal G}\otimes_{_{{\cal E}}}{\cal E}(\Lambda ^{k+1} _{_{\R}}T^*_X)
$$
par la formule classique
\begin{eqnarray*} 
&\displaystyle{\nabla_{_{\cal G}}\omega\, (\xi _0,...,\xi _k):=\sum_{0\leq j \leq k}(-1)^j\nabla_{_{\cal G}}(\omega (\xi _0,...,\widehat{\xi _j},..., \xi _k))(\xi _j)+ }&
\\
\\
&\displaystyle{+\sum_{0\leq j<l \leq k}(-1)^{j+l}\omega ([\xi_j,\xi_l],\xi _0,...,\widehat{\xi _j},...,\widehat{\xi _l},..., \xi _k)   }&
\end{eqnarray*}
pour tout $\omega \in({\cal G}\otimes_{_{{\cal E}}}{\cal E}(\Lambda ^k_{_{\R}}T^*_X))(U)$ et tout champ de vecteurs complexes $\xi _0,...,\xi _k\in{\cal E}(T_X\otimes_{_{\R}}\C)(U)$
L'extension ainsi définie vérifie la règle de Leibnitz $\nabla_{_{\cal G}}(g\otimes f)=\nabla_{_{\cal G}}\,g\wedge f+g\otimes df$ pour tout $g\in{\cal G}(U)$ et $f\in{\cal E}(\Lambda ^k_{_{\R}} T^*_X)(U)$. En effet
\begin{eqnarray*} 
&\displaystyle{\nabla_{_{\cal G}}(g\otimes f)\, (\xi _0,...,\xi _k):=
\sum_{0\leq j \leq k}(-1)^j\nabla_{_{\cal G}}(g\cdot f (\xi _0,...,\widehat{\xi _j},..., \xi _k))(\xi _j)+ }&
\\
\\
&\displaystyle{+\sum_{0\leq j<l \leq k}(-1)^{j+l}g\cdot f([\xi_j,\xi_l],\xi _0,...,\widehat{\xi _j},...,\widehat{\xi _l},..., \xi _k)  = }&
\\
\\
&\displaystyle{=
\sum_{0\leq j \leq k}(-1)^j\Big[\nabla_{_{\cal G}}g(\xi _j)\cdot f (\xi _0,...,\widehat{\xi _j},..., \xi _k)+g\cdot (\xi _j\,.f(\xi _0,...,\widehat{\xi _j},..., \xi _k))\Big]+}&
\\
\\
&\displaystyle{+\sum_{0\leq j<l \leq k}(-1)^{j+l}g\cdot f([\xi_j,\xi_l],\xi _0,...,\widehat{\xi _j},...,\widehat{\xi _l},..., \xi _k)  = }&
\\
\\
&\displaystyle{=(\nabla_{_{\cal G}}\,g \wedge f+g\otimes df)(\xi _0,...,\xi _k)}.&
\end{eqnarray*}
Le fait que $d^2=0$ entraîne l'existence du tenseur de courbure de $\nabla_{{\cal G}}$
$$
\Theta(\nabla_{_{\cal G}})\in \Big({\cal E}nd_{_{{\cal E}(\C)}}({\cal G}) \otimes_{_{{\cal E}(\C) }}{\cal E}(\Lambda ^2_{_{\R}} T^*_X)\Big)(X)
$$
définie par la formule $\Theta(\nabla_{_{\cal G}})(\xi ,\eta)\cdot g:=(\nabla^2_{_{\cal G}}\,g)(\xi ,\eta)$ pour tout $\xi ,\eta\in {\cal E}(T_X)(U)$ et $g\in {\cal G}(U)$. On note de plus par $\xi _{_{\nabla_{_{\cal G}}} }\,.g:=\nabla_{_{\cal G}}g(\xi)$ la dérivée covariante de la section $g$ le long du champ de vecteurs $\xi$. La définition de l'extension de la connexion $\nabla_{_{\cal G}}$ implique de façon immédiate la formule
$$
\xi _{_{\nabla_{_{\cal G}}} }\,.\,(\eta_{_{\nabla_{_{\cal G}}} }\,.\,g)-\eta_{_{\nabla_{_{\cal G}}} }\,.\,(\xi _{_{\nabla_{_{\cal G}}} }\,.\,g)=[\xi ,\eta]_{_{\nabla_{_{\cal G}}} }\,.\,g+\Theta(\nabla_{_{\cal G}})(\xi ,\eta)\cdot g.
$$
Le tenseur de courbure $\Theta(\nabla_{_{\cal G}})$ de la connexion $\nabla_{_{\cal G}}$ mesure donc le défaut de commutation des dérivées covariantes secondes des sections de ${\cal G}$. Il est aussi élémentaire de vérifier l'identité
\begin{eqnarray}\label{formulcurvat} 
\nabla^2_{_{\cal G}}\,\omega =\Theta(\nabla_{_{\cal G}})\wedge \omega 
\end{eqnarray}
pour tout $\omega \in({\cal G}\otimes_{_{{\cal E}}}{\cal E}(\Lambda ^{\bullet}_{_{\R}}T^*_X))(U)$.
Le fait que les opérateurs $\theta_{_J }$ et $\bar{\theta }_{_J }$ vérifient la règle de Leibnitz entraîne que 
$$
\theta_{_J }\in {\cal H} om_{_{{\cal E}(\C) } }({\cal E}^{p,q}_{_{X,J} }  , {\cal E}^{p+2,q-1}_{_{X,J} } )(X)\quad\mbox{et} \quad
\bar{\theta }_{_J }\in {\cal H} om_{_{{\cal E}(\C) } }({\cal E}^{p,q}_{_{X,J} }  , {\cal E}^{p-1,q+2}_{_{X,J} } )(X)
$$
On définit alors les opérateurs de torsion sur ${\cal G}$
\begin{eqnarray*}
\theta _{_{{\cal G},J} } :=
\I_{_{\cal G} } \otimes_{_{{\cal E}(\C)}}\theta_{_J }:{\cal G}\otimes_{_{{\cal E}(\C)}}{\cal E}^{p,q}_{_{X,J} } 
\longrightarrow{\cal G}\otimes_{_{{\cal E}(\C)}}{\cal E}^{p+2,q-1}_{_{X,J} } 
\\
\\
\bar{\theta }_{_{{\cal G},J} } :=
\I_{_{\cal G} } \otimes_{_{{\cal E}(\C)}}\bar{\theta }_{_J }:{\cal G}\otimes_{_{{\cal E}(\C)}}{\cal E}^{p,q}_{_{X,J} } 
\longrightarrow{\cal G}\otimes_{_{{\cal E}(\C)}}{\cal E}^{p-1,q+2}_{_{X,J} } .
\end{eqnarray*}
De façon explicite ces opérateurs sont définis de façon analogue aux opérateurs $\theta_{_J }$ et $\bar{\theta }_{_J }$. Ce sont des dérivations, autrement dit on a les formules
\begin{eqnarray*}
\theta _{_{{\cal G},J} }(\omega \wedge f)=\theta _{_{{\cal G},J} }\,\omega \wedge f+(-1)^{\deg\, \omega }\,\omega \wedge  \theta_{_J }f
\\
\\
\bar{\theta }_{_{{\cal G},J} } (\omega \wedge f)=\bar{\theta }_{_{{\cal G},J} } \,\omega \wedge f+(-1)^{\deg\, \omega }\,\omega \wedge \bar{\theta }_{_J } f
\end{eqnarray*}
pour tout $\omega \in ({\cal G}\otimes_{_{{\cal E}}}{\cal E}(\Lambda ^{\bullet}_{_{\R}} T^*_X)(U)$ et $f\in{\cal E}(\Lambda ^{\bullet}_{_{\R}} T^*_X)(U)$. Comme dans le cas de la différentielle extérieure on a la décomposition 
$$
\nabla_{_{{\cal G}}}=\nabla^{1,0}_{_{{\cal G},J} }\,+\,\nabla^{0,1}_{_{{\cal G},J} }\,-\,\theta _{_{{\cal G},J} }\,-\,
\bar{\theta }_{_{{\cal G},J} }
$$
où les opérateurs
$$
\nabla^{1,0}_{_{{\cal G},J} } :{\cal G}\otimes_{_{{\cal E}(\C)}}{\cal E}^{p,q}_{_{X,J} } 
\longrightarrow{\cal G}\otimes_{_{{\cal E}(\C)}}{\cal E}^{p+1,q}_{_{X,J} } 
$$
et
$$
\nabla^{0,1}_{_{{\cal G},J} } :{\cal G}\otimes_{_{{\cal E}(\C)}}{\cal E}^{p,q}_{_{X,J} } 
\longrightarrow{\cal G}\otimes_{_{{\cal E}(\C)}}{\cal E}^{p,q+1}_{_{X,J} } 
$$
sont définis par les formules analogues à celles qui définisent les opérateurs $\bar{\partial}_{_J } $ et $\partial_{_J }$,
\begin{eqnarray} 
&\displaystyle{\nabla^{1,0}_{_{{\cal G},J} }\,\omega\, (\xi _0,...,\xi _k):=
\sum_{0\leq j \leq k}(-1)^j \nabla_{_{{\cal G}} }(\omega (\xi _0,...,\widehat{\xi _j},..., \xi _k))(\xi^{1,0}_j)+ }&\nonumber
\\\nonumber
\\
&\displaystyle{+\sum_{0\leq j<l \leq k}(-1)^{j+l}\omega ([\xi^{1,0}_j,\xi^{1,0}_l]^{1,0}+[\xi^{0,1} _j,\xi^{1,0} _l]^{0,1}+[\xi^{1,0} _j,\xi^{0,1} _l]^{0,1},\xi _0,...,\widehat{\xi _j},...,\widehat{\xi _l},..., \xi _k)   }&\label{extendcon1,0} 
\end{eqnarray}
et
\begin{eqnarray} 
&\displaystyle{\nabla^{0,1}_{_{{\cal G},J} }\,\omega\, (\xi _0,...,\xi _k):=
\sum_{0\leq j \leq k}(-1)^j\nabla_{_{{\cal G}} }(\omega (\xi _0,...,\widehat{\xi _j},..., \xi _k))(\xi^{0,1}_j) + }&\nonumber
\\\nonumber
\\
&\displaystyle{+\sum_{0\leq j<l \leq k}(-1)^{j+l}\omega ([\xi^{0,1}_j,\xi^{0,1}_l]^{0,1}+[\xi^{0,1} _j,\xi^{1,0} _l]^{1,0}+[\xi^{1,0} _j,\xi^{0,1} _l]^{1,0},\xi _0,...,\widehat{\xi _j},...,\widehat{\xi _l},..., \xi _k)   }&\label{extendcon0,1}
\end{eqnarray}
Le fait que $\nabla_{_{{\cal G}}}$ vérifie la règle de Leibnitz implique les formules
\begin{eqnarray*}
\nabla^{1,0}_{_{{\cal G},J}}\,(g\otimes f)=\nabla^{1,0}_{_{{\cal G},J}}\,g\wedge f+g\otimes \partial_{_J } f
\\
\\
\nabla^{0,1}_{_{{\cal G},J}}\,(g\otimes f)=\nabla^{0,1}_{_{{\cal G},J}}\,g\wedge f+g\otimes \bar{\partial}_{_J } f
\end{eqnarray*}
pour tout $g\in {\cal G}(U)$ et $f\in{\cal E}(\Lambda ^{\bullet}_{_{\R}} T^*_X)(U)$. En degré zéro on a les formules
$$
\nabla^{1,0}_{_{{\cal G},J}}\,g=\frac{1}{2}\Big(\nabla_{_{{\cal G}}}\,g-i(\nabla_{_{{\cal G}}}\,g)\circ J\Big)\quad\mbox{et}\quad
\nabla^{0,1}_{_{{\cal G},J}}\,g=\frac{1}{2}\Big(\nabla_{_{{\cal G}}}\,g+i(\nabla_{_{{\cal G}}}\,g)\circ J\Big)  
$$
pour tout $g\in {\cal G}(U)$. En général on a la définition suivante.
\begin{defi} 
Soit ${\cal G}$ un faisceau de ${\cal E}(\C)$-modules sur $X$. Une connexion de type $(0,1)$ sur le faisceau ${\cal G}$ est un morphisme de faisceaux de groupes additifs $\nabla''_{_{\cal G}}:{\cal G}\longrightarrow{\cal G}\otimes_{_{{\cal E}(\C)}}{\cal E}^{0,1}_{_{X,J} } $ tel que 
$\nabla''_{_{\cal G}}(g\cdot f)=\nabla''_{{\cal G}}\,g\cdot f+g\otimes \bar{\partial}_{_J }f$ pour tout $g\in{\cal G}(U)$ et $f\in {\cal E}(\C)(U)$, où $U\subset X$ est un ouvert quelconque.
\end{defi} 
On a bien sûr une définition analogue pour les connexions de type $(1,0)$. Comme précédemment une connexion de type $(0,1)$, (resp. $(1,0)$) peut être étendue grâce à la formule $\ref{extendcon0,1}$, (resp. $\ref{extendcon1,0}$) ou grâce à la règle de Leibnitz.
On rappelle maintenant que si $A$ et $B$ sont deux endomorphismes du faisceau de ${\cal E}(\C)$-modules ${\cal G}\otimes_{_{{\cal E}}}{\cal E}(\Lambda ^{\bullet} _{_{\R}} T^*_X)$, leur crochet de commutation est défini par la formule $[A,B]:=AB-(-1)^{\deg\,A\cdot \deg\, B}BA$. La décomposition précédente de $\nabla_{_{{\cal G}}}$ implique la décomposition suivante au niveau des opérateurs,
\begin{eqnarray*}
&\displaystyle{\nabla^2_{_{{\cal G}}}=(\nabla^{1,0}_{_{{\cal G},J} }\,+\,\nabla^{0,1}_{_{{\cal G},J} }\,-\,\theta _{_{{\cal G},J} }\,-\,
\bar{\theta }_{_{{\cal G},J} })^2=} &
\end{eqnarray*}
\begin{eqnarray*}
&\displaystyle{=\underbrace{(\nabla^{1,0}_{_{{\cal G},J} })^2-[\nabla^{0,1}_{_{{\cal G},J} }\,,\theta _{_{{\cal G},J} }]}_{2,0}+
\underbrace{(\nabla^{0,1}_{_{{\cal G},J} })^2-[\nabla^{1,0}_{_{{\cal G},J} }\,,\bar{\theta } _{_{{\cal G},J} }]}_{0,2} +
\underbrace{\theta ^2_{_{{\cal G},J} }}_{4,-2} +\underbrace{\bar{\theta }^2_{_{{\cal G},J} }}_{-2,4}+  } & \\
\\
&\displaystyle{+\underbrace{[\nabla^{1,0}_{_{{\cal G},J} }\,,\nabla^{0,1}_{_{{\cal G},J} }]
+[\theta _{_{{\cal G},J}}\,,\bar{\theta } _{_{{\cal G},J} }]}_{1,1}
- \underbrace{[\nabla^{1,0}_{_{{\cal G},J} }\,,\theta _{_{{\cal G},J} }]}_{3,-1} 
 - \underbrace{[\nabla^{0,1}_{_{{\cal G},J} }\,,\bar{\theta } _{_{{\cal G},J} }]}_{-1,3}}.& 
\end{eqnarray*}
D'autre part en considérant la décomposition de la forme de courbure 
$$
\Theta(\nabla_{_{\cal G}})=\Theta(\nabla_{_{\cal G}})^{2,0}_{_{J} }+  \Theta(\nabla_{_{\cal G}})^{1,1}_{_{J} }+  \Theta(\nabla_{_{\cal G}})^{0,2}_{_{J} }
$$
en ses composantes de type $(2,0),\,(1,1),\,(0,2)$ et la formule $\eqref{formulcurvat}$ on déduit les identités suivantes au sens des opérateurs
\begin{align*}
\Theta(\nabla_{_{\cal G}})^{2,0}_{_{J} }\wedge \cdot=(\nabla^{1,0}_{_{{\cal G},J} })^2-[\nabla^{0,1}_{_{{\cal G},J} }\,,\theta _{_{{\cal G},J} }]                                                 &\qquad \qquad  
\Theta(\nabla_{_{\cal G}})^{0,2}_{_{J} }\wedge \cdot=(\nabla^{0,1}_{_{{\cal G},J} })^2-[\nabla^{1,0}_{_{{\cal G},J} }\,,\bar{\theta } _{_{{\cal G},J} }]
\\
\\
\Theta(\nabla_{_{\cal G}})^{1,1}_{_{J} }\wedge \cdot=[\nabla^{1,0}_{_{{\cal G},J} }\,,\nabla^{0,1}_{_{{\cal G},J} }]
+[\theta _{_{{\cal G},J}}\,,\bar{\theta } _{_{{\cal G},J} }]           &\qquad \qquad
\theta ^2_{_{{\cal G},J} }=0,\quad \bar{\theta }^2_{_{{\cal G},J} }=0
\\
\\
[\nabla^{1,0}_{_{{\cal G},J} }\,,\theta _{_{{\cal G},J} }]=0           &\qquad \qquad
[\nabla^{0,1}_{_{{\cal G},J} }\,,\bar{\theta } _{_{{\cal G},J} }]=0.
\end{align*}
En particulier si ${\cal G}={\cal E}(\C)$ et $\nabla_{_{\cal G}}=d$ on a les identités supplémentaires 
$$
\partial^2_{_J }=[\bar{\partial}_{_J },\theta_{_J }],\qquad \bar{\partial}^2_{_J }=[\partial_{_J },\bar{\theta }_{_J }]\quad\mbox{et}\quad [\partial_{_J },\bar{\partial}_{_J }]=-[\theta_{_J },\bar{\theta }_{_J }].
$$
En conclusion on a les identités fondamentales de la géométrie presque complexe: 
\begin{eqnarray*}
&\displaystyle{\partial^2_{_J }=
\bar{\partial}_{_J }\theta_{_J }+\theta_{_J }\bar{\partial}_{_J },\qquad 
\bar{\partial}^2_{_J }=
\partial_{_J }\bar{\theta }_{_J }+\bar{\theta }_{_J }\partial_{_J },}&
\\
\\
&\displaystyle{\partial_{_J }\bar{\partial}_{_J }+\bar{\partial}_{_J }\partial_{_J }=
-\theta_{_J }\bar{\theta }_{_J }-\bar{\theta }_{_J }\theta_{_J },}&
\\
\\
&\displaystyle{\partial_{_J }\theta_{_J }=-\theta_{_J }\partial_{_J },\qquad 
\bar{\partial}_{_J }\bar{\theta }_{_J }=-\bar{\theta }_{_J }\bar{\partial}_{_J },}&
\\
\\
&\displaystyle{\theta^2_{_J }=0,\qquad\bar{\theta }^2_{_J }=0.}&
\end{eqnarray*}  
En général en degré zéro on a les formules
\begin{eqnarray}
&\displaystyle{
\Theta(\nabla_{_{\cal G}})^{2,0}_{_{J} }=
(\nabla^{1,0}_{_{{\cal G},J} })^2-\theta _{_{{\cal G},J} }\,\nabla^{0,1}_{_{{\cal G},J} },  }&\label{2,0} 
\\\nonumber
\\
&\displaystyle{\Theta(\nabla_{_{\cal G}})^{0,2}_{_{J} }=
(\nabla^{0,1}_{_{{\cal G},J} })^2-\bar{\theta }_{_{{\cal G},J} }\,\nabla^{1,0}_{_{{\cal G},J} },
}& \label{0,2} 
\\\nonumber
\\
&\displaystyle{
\Theta(\nabla_{_{\cal G}})^{1,1}_{_{J} }=[\nabla^{1,0}_{_{{\cal G},J} }\,,\nabla^{0,1}_{_{{\cal G},J} }],
}&\label{1,1} 
\end{eqnarray}
qui sont équivalentes aux identités évidentes
\begin{eqnarray*}
\xi^{1,0}_{_{\nabla_{_{\cal G}}} }\,.\,(\eta^{1,0}_{_{\nabla_{_{\cal G}}} }\,.\,g)-\eta^{1,0}_{_{\nabla_{_{\cal G}}} }\,.\,(\xi^{1,0} _{_{\nabla_{_{\cal G}}} }\,.\,g)=[\xi^{1,0} ,\eta^{1,0}]_{_{\nabla_{_{\cal G}}} }\,.\,g+\Theta(\nabla_{_{\cal G}})^{2,0}_{_{J}}(\xi^{1,0} ,\eta^{1,0})\cdot g,
\\
\\
\xi^{0,1}_{_{\nabla_{_{\cal G}}} }\,.\,(\eta^{0,1}_{_{\nabla_{_{\cal G}}} }\,.\,g)-\eta^{0,1}_{_{\nabla_{_{\cal G}}} }\,.\,(\xi^{0,1} _{_{\nabla_{_{\cal G}}} }\,.\,g)=[\xi^{0,1} ,\eta^{0,1}]_{_{\nabla_{_{\cal G}}} }\,.\,g+\Theta(\nabla_{_{\cal G}})^{0,2}_{_{J}}(\xi^{0,1} ,\eta^{0,1})\cdot g,
\\
\\
\xi^{1,0}_{_{\nabla_{_{\cal G}}} }\,.\,(\eta^{0,1}_{_{\nabla_{_{\cal G}}} }\,.\,g)-
\eta^{0,1}_{_{\nabla_{_{\cal G}}} }\,.\,(\xi^{1,0} _{_{\nabla_{_{\cal G}}} }\,.\,g)=
[\xi^{1,0} ,\eta^{0,1}]_{_{\nabla_{_{\cal G}}} }\,.\,g+\Theta(\nabla_{_{\cal G}})^{1,1}_{_{J}}(\xi^{1,0} ,\eta^{0,1})\cdot g.
\end{eqnarray*}
On a donc en particulier que la composante de type $(2,0)$, (resp. $(0,2)$) du tenseur de courbure mesure le défaut de commutation des dérivées covariantes secondes des sections de ${\cal G}$ le long des champs de vecteurs de type $(1,0)$, (resp. $(0,1)$). La composante de type $(1,1)$ du tenseur de courbure exprime  le défaut de commutation des dérivées covariantes secondes des sections de ${\cal G}$ le long des champs de vecteurs de type $(1,0)$ et  $(0,1)$. Soit ${\cal G}$ un faisceau de ${\cal E}(\C)$-modules localement de type fini, soit $\psi\equiv (\psi_1,...,\psi_r)\in {\cal G}^{\oplus r}(U)$ un système de générateurs locaux et 
$$
\omega =\psi\cdot f \in({\cal G}\otimes_{_{{\cal E}}}{\cal E}(\Lambda ^k_{_{\R}}T^*_X))(U),\qquad f\in M_{r,1}({\cal E}(\Lambda ^k_{_{\R}}T^*_X)(U)).
$$
Soient de plus $A\in M_{r,r} ({\cal E}(T^*_X)(U)),\,A'_{_{J}}\in M_{r,r}({\cal E}^{1,0}_{_{X,J} } (U)),\,A''_{_{J}}\in M_{r,r}({\cal E}^{0,1}_{_{X,J} } (U))$ telles que $\nabla_{_{\cal G}}\,\psi =\psi\cdot A$ et $A=A'_{_{J}}+A''_{_{J}}$. La règle de Leibnitz implique alors les égalités
$$
\nabla_{_{\cal G}}\,\omega =\psi\cdot (df+A\wedge f)\quad\mbox{et}\quad  \Theta(\nabla_{_{\cal G}})\wedge \omega =\psi\cdot (dA+A\wedge A)\wedge f.
$$
De plus on a les identités 
\begin{align*}
\nabla^{1,0}_{_{{\cal G},J} }\,\omega =\psi\cdot(\partial_{_J }f+A'_{_{J}}\wedge f),          &\qquad \qquad   
\nabla^{0,1}_{_{{\cal G},J} }\,\omega =\psi\cdot (\bar{\partial}_{_J }f+A''_{_{J}}\wedge f),
\\
\\
\theta _{_{{\cal G},J} }\,\omega =\psi\cdot\theta_{_J }\,f,                        &\qquad \qquad 
\bar{\theta }_{_{{\cal G},J} }\,\omega =\psi\cdot \bar{\theta }_{_J }\,f.
\end{align*}
En décomposant la $2$-forme $dA+A\wedge A$ où en explicitant les identités $\eqref{2,0},\,\eqref{0,2}$ et $\eqref{1,1}$ on obtient les expressions locales  suivantes.
\begin{eqnarray*}
&\displaystyle{\Theta(\nabla_{_{\cal G}})^{2,0}_{_{J} }\wedge \omega =
\psi\cdot (\partial_{_J }A'_{_{J}}+A'_{_{J}}\wedge A'_{_{J}}-\theta_{_J }\,A''_{_{J}})\wedge f}&
\\
\\
&\displaystyle{\Theta(\nabla_{_{\cal G}})^{0,2}_{_{J} }\wedge \omega =
\psi\cdot (\bar{\partial}_{_J }A''_{_{J}}+A''_{_{J}}\wedge A''_{_{J}}-\bar{\theta }_{_J }\,A'_{_{J}})\wedge f}&
\\
\\
&\displaystyle{\Theta(\nabla_{_{\cal G}})^{1,1}_{_{J} }\wedge \omega =
\psi\cdot (\bar{\partial}_{_J }A'_{_{J}}+\partial_{_J }A''_{_{J}}+A'_{_{J}}\wedge A''_{_{J}}+A''_{_{J}}\wedge A'_{_{J}})\wedge f,}&
\end{eqnarray*}
\section{Connexions hermitiennes sur les fibrés vectoriels au dessus des variétés presque complexes}\label{basic-herm-con}
Nous considérons à partir de maintenant un fibré vectoriel complexe $\ci$, $F\longrightarrow X$ et ${\cal G}={\cal E}(F):=$faisceau des sections $\ci$ de $F$. Soit $h\in{\cal E}(F^*\otimes_{_{\C}}\overline{F}^* )(X)$ une métrique hermitienne sur $F$. On rappelle qu'une connexion 
$$
\nabla_{_F}:{\cal E}(F)\longrightarrow{\cal E}(F)\otimes_{_{{\cal E}}}{\cal E}(T^*_X)\simeq{\cal E}(F)\otimes_{_{{\cal E}(\C)}}{\cal E}(T^*_X\otimes_{_{\R}}\C)
$$
sur $F$ est dite $h$-hermitienne si pour tout champ de vecteurs complexes 
$\xi \in{\cal E}(T_X\otimes_{_{\R}}\C)(U)$ et toute sections $\sigma ,\tau\in{\cal E}(F)(U)$, ($U\subseteq X$ est un ouvert quelconque), on a la formule
$$
\xi .h(\sigma ,\tau)=h(\xi _{_{\nabla} }.\sigma ,\tau)+h(\sigma ,\bar{\xi} _{_{\nabla} }.\tau).
$$
Il est bien sûr équivalent de restreindre l'identité précédente aux seuls champs de vecteurs 
$\xi\in{\cal E}(T^{1,0}_{X,J})(U)$. On a alors que la donnée d'une connexion 
$$
\nabla''_{_F}:{\cal E}(F)\longrightarrow{\cal E}(F)\otimes_{_{{\cal E}(\C)}}{\cal E}^{0,1}_{_{X,J} } 
$$ 
de type $(0,1)$ entraîne l'existence d'une unique connexion $h$-hermitienne $\nabla_{_F}$ sur le fibré $F$ telle que $\nabla^{0,1} _{_F}=\nabla''_{_F}$. En effet la partie de type $(1,0)$ de $\nabla_{_F}$ est donnée par la formule
$$
h(\nabla^{1,0} _{_F}\sigma(\xi ) ,\tau)=\xi .h(\sigma ,\tau)-h(\sigma ,\nabla''_{_F}\tau(\bar{\xi }) )
$$
pour tout $(1,0)$-champ de vecteurs $\xi\in{\cal E}(T^{1,0}_{X,J})(U)$ et toutes sections $\sigma ,\tau\in{\cal E}(F)(U)$.
Bien évidemment on a un résultat analogue pour les connexions de type $(1,0)$. Soit $(e_1,...,e_r)\in {\cal E}(F)^{\oplus r}(U)$ un repère de $F_{|_U}$. On a l'identification $\nabla_{_F}\simeq_{e}d+A $ par rapport au repère $(e_1,...,e_r)$. Soit de plus $H:=(h(e_{\lambda},e_{\mu}))_{\lambda ,\mu }$ la matrice hermitienne de la métrique $h$. Le fait que la connexion $\nabla_{_F}$ soit $h$-hermitienne  équivaut localement  aux égalités
$$
\xi .H_{\lambda ,\mu }=\sum_{1\leq s\leq r}\Big(A' _{s,\lambda }(\xi )H_{s,\mu}+\overline{A''_{s,\mu}(\bar{\xi} )}H_{\lambda ,s}\Big),
$$
$\xi\in{\cal E}(T^{1,0}_{X,J})(U)$. On a alors avec des notations matricielles la relation $\partial_{_J }H=A_{_{J} } '^{t}H+H\overline{A''_{_{J} } }$. Le fait que la matrice $H$ soit hermitienne implique que cette relation est équivalente à la relation
\begin{eqnarray}\label{hermitiancon} 
A'_{_{J}}=\overline{H}^{-1}(\partial_{_J }\overline{H}-\overline{A''_{_{J} }}^{t}\,\overline{H} ).  
\end{eqnarray}
On a en conclusion qu'une connexion $\nabla_{_F}$ est $h$-hermitienne si et seulement si la relation $\eqref{hermitiancon}$ est satisfaite sur tout les ouverts de trivialisation de $F$. Considérons maintenant le produit sesquilinéaire
$$
\{\cdot,\cdot\}_h:  {\cal E}(\Lambda ^p_{_{\R}} T^*_X\otimes_{_{\R}} F)\times {\cal E}(\Lambda ^q_{_{\R}} T^*_X\otimes_{_{\R}} F)\longrightarrow {\cal E}(\Lambda ^{p+q} _{_{\R}} T^*_X\otimes_{_{\R}} \C)
$$
sur le faisceau ${\cal E}(\Lambda ^{\bullet} _{_{\R}} T^*_X\otimes_{_{\R}} F)$ défini par la formule
$$
\{\sigma ,\tau\}_h(\xi )=\sum_{|I|=p}\varepsilon (I)h(\sigma (\xi _{I} ),\tau(\bar{\xi }_{\complement I}  )),  
$$
où $\xi =(\xi _1,...,\xi _{p+q}), \,\xi _j\in{\cal E}(T_X\otimes_{_{\R}}\C)(U)$ et $\varepsilon (I)$ désigne le signe de la permutation $(1,...,p+q)\rightarrow (I,\complement I)$. Alors le fait que la connexion $\nabla_{_F}$ soit hermitienne est équivalent à l'identité plus générale 
$$
d\{\sigma ,\tau\}_h=\{\nabla_{_F}\sigma ,\tau\}_h+(-1)^{\deg\,\sigma } \{\sigma ,\nabla_{_F}\tau\}_h
$$ 
qui équivaut aussi à une des identités
\begin{eqnarray*}
\partial_{_J }\{\sigma ,\tau\}_h=\{\nabla^{1,0}_{_{F,J} }\sigma ,\tau\}_h+
(-1)^{\deg\,\sigma } \{\sigma ,\nabla^{0,1}_{_{F,J} }\tau\}_h
\\
\\
\bar{\partial}_{_J }\{\sigma ,\tau\}_h=\{\nabla^{0,1}_{_{F,J} }\sigma ,\tau\}_h+
(-1)^{\deg\,\sigma } \{\sigma ,\nabla^{1,0}_{_{F,J} }\tau\}_h.
\end{eqnarray*}
On obtient alors, en appliquant la différentielle extérieure à la première des trois identités précédentes, l'identité 
$0=\{\Theta (\nabla_{_F})\sigma ,\tau\}_h+\{\sigma ,\Theta (\nabla_{_F})\tau\}_h$ qui implique, pour des raisons de bidegré, l'identité 
$$
0=\{\Theta (\nabla_{_F})^{1,1}_{_{J} }  \sigma ,\tau\}_h+\{\sigma ,\Theta (\nabla_{_F})^{1,1}_{_{J} }\tau\}_h.
$$
Si $\deg \,\sigma =\deg\,\tau=0$ on déduit l'égalité 
\begin{eqnarray}\label{hermitcurvat} 
0=h\Big(\Theta (\nabla_{_F})^{1,1}_{_{J} } (\xi ,\eta)\cdot \sigma ,\tau\Big)+h\Big(\sigma ,\Theta (\nabla_{_F})^{1,1}_{_{J} } (\bar{\xi }  ,\bar{\eta} )\cdot \tau\Big)
\end{eqnarray}
qui montre que pour tout champ de vecteurs réels $\xi ,\,\eta\in {\cal E}(T_X)(U)$ on a 
$$
i\Theta (\nabla_{_F})^{1,1}_{_{J} }(\xi ,\eta) \in{\cal E}(\mbox{Herm}_{h}(F))(U),
$$
où $\mbox{Herm}_{h}(F)$ désigne le fibré (réel) des endomorphismes $h$-hermitiens de $F$. 
Considérons maintenant l'expression locale de la composante de type $(1,1)$ du tenseur de courbure 
$$
\Theta (\nabla_{_F})^{1,1}_{_{J} } =\sum_{1\leq \lambda ,\mu \leq r}C_{\lambda ,\mu }\otimes e^*_{\mu }\otimes e_{\lambda }
$$
de la connexion hermitienne $\nabla_{_F}$. On a 
$$
C:=\bar{\partial}_{_J }A'_{_{J}}+\partial_{_J }A''_{_{J}}+A'_{_{J}}\wedge A''_{_{J}}+A''_{_{J}}\wedge A'_{_{J}}.
$$
Si 
$(\zeta _k)_k\in {\cal E}(T^{1,0}_{X,J})^{\oplus n}(U)$ est un repère du fibré $T^{1,0}_{U,J}$, on a l'expression locale suivante 
$$
\Theta (\nabla_{_F})^{1,1}_{_{J} } =\sum_{\substack{1\leq \lambda ,\mu \leq r
\\
\\
1\leq k,l\leq n }} C^{k,l}_{\lambda ,\mu}\,\zeta ^*_k\wedge \bar{\zeta } ^*_l\otimes  e^*_{\mu }\otimes e_{\lambda }.   
$$
L'identité $\eqref{hermitcurvat}$ entraîne que
si en un point $x_0\in U$ le repère $e_1(x_0),...,e_r(x_0)$ est $h(x_0)$-orthonormé alors on a les relations $\overline{C^{k,l}_{\lambda ,\mu}}(x_0)=C^{l,k}_{\mu ,\lambda }(x_0)$. Si de plus 
$\nabla_{_F}^{0,1}e_k(x_0)=0 $ pour tout $k$, on obtient en utilisant l'expression $\eqref{hermitiancon}$ l'égalité 
\begin{eqnarray}\label{ponctcurvat}
C(x_0)=(\bar{\partial}_{_J }\partial_{_J }\overline{H}-\bar{\partial}_{_J }\overline{H}\wedge \partial_{_J }\overline{H} 
+\partial_{_J }A''_{_{J}}-\bar{\partial}_{_J }\overline{A''_{_{J}}}^t)(x_0).
\end{eqnarray}
\section{Extension de l'opérateur $\bar{\partial}_{_J }$ aux puissances de Schur du fibré des $(1,0)$-formes}
Rappelons qu'étant donné un espace vectoriel complexe $V$ de dimension complexe $r$, les représentations irréductibles de $GL_{_\C}(V)$ sont en correspondance biunivoque avec le plus haut poids $\lambda =(\lambda _1, ...\lambda _r),\,\lambda _1\geq \lambda _2\geq ...\geq \lambda _r$ de la représentation d'un sous-tore maximal $T^r\simeq (\C^*)^r<GL_{_\C}(V)$, $(t_1,...,t_r)\mapsto t_1^{\lambda _1}\cdot \cdot\cdot t_r^{\lambda _r}$. On note $S^{\lambda}V$ l'espace de la représentation associée, qu'on appelle puissance de Schur associée au poids $\lambda$. On a par exemple 
\begin{eqnarray*}
&\displaystyle{
S^{(m,0,...,0)}V=S^mV }&\mbox{puissance symétrique usuelle}
\\
&\displaystyle{
S^{(1,1,...,1,0,...,0)}V=\Lambda ^kV  }&  \mbox{puissance extérieure}. 
\end{eqnarray*}  
Nous renvoyons le lecteur aux ouvrages classiques de \cite{Fu-Ha} pour une explication détaillé de la notion de puissance de Schur.
\\
\\
Considérons maintenant les connexions de type $(0,1)$ 
$$
\bar{\partial}_{_{J,p} } :=(-1)^p\bar{\partial}_{_{J}}
:{\cal E}^{p,0}_{_{X,J} }
\longrightarrow {\cal E}^{p,0}_{_{X,J} }\otimes_{_{{\cal E}(\C)}}{\cal E}^{0,1}_{_{X,J} } 
$$
sur les fibrés $\Lambda ^{p,0}_{_J}T_X^*$. De façon explicite les connections $\bar{\partial}_{_{J,p} }$ sont définies par les formules
\begin{eqnarray}\label{(0,1)conn-sur-(p,0)} 
&\displaystyle{
\left<\bar{\partial}_{_{J,p} }\omega (\eta), \xi _1,...,\xi _p\right>:=\bar{\partial}_{_{J} }\omega (\eta, \xi _1,...,\xi _p)}&\nonumber
\\\nonumber
\\
&\displaystyle{=\eta\,.\,\omega (\xi _1,...,\xi _p)+\sum_{1\leq l \leq p}(-1)^l\omega ([\eta, \xi _l]^{1,0}, \xi _1,...,\widehat{\xi _l},..., \xi _p)}& 
\end{eqnarray}
pour tout $\omega \in {\cal E}^{p,0}_{_{X,J} }(U),\,\eta\in{\cal E}(T^{0,1}_{X,J})(U)$ et 
$\xi _1,...,\xi _p\in{\cal E}(T^{1,0}_{X,J})(U)$. Bien évidement dans le cas complexe intégrable les faisceaux de ${\cal O}_{X,J}$-modules 
$
\Omega^p_{X,J}\equiv {\cal O}(\Lambda ^{p,0}_{_J}T_X^*):=\Ker\bar{\partial}_{_{J,p} }
$
sont localement libres et donnent une structure de fibré vectoriel holomorphe aux fibrés  $\Lambda ^{p,0}_{_J}T_X^*$. De plus on a l'identité
$$
\bar{\partial}_{_{J,p}}=\I_{_{\Omega^p _{X,J}}} \otimes_{_{{\cal O}_X } } \bar{\partial}_{_{J} }.
$$
Dans le cas presque complexe, en étendant la connexion $\bar{\partial}_{_{J,1}}$  à toutes les puissances de Schur $F^{\lambda }_{_{J}}:=S^{\lambda }\Lambda ^{1,0}_{_J}T_X^* $ on obtient des connections de type $(0,1)$ canoniques 
$$
\bar{\partial}_{{F^{\lambda }_{_{J}} } }: {\cal E}(F^{\lambda  }_{_{J}}) \longrightarrow {\cal E}(\Lambda ^{0,1}_{_J}T^*_{X}\otimes_{_{\C}}F^{\lambda  }_{_{J}})
$$
sur les fibrés $F^{\lambda }_{_{J}}$.
De façon analogue, la connexion  induite sur $T^*_{X,J}$ par $\bar{\partial}_{_{J,1}}$ grâce au $\C$-isomorphisme canonique de $\Lambda ^{1,0}_{_J}T_X^* $ avec $T^*_{X,J}$ peut être étendue aux puissances de Schur $S^{\lambda }T^*_{X,J}$. Pour simplifier nous désignerons aussi $S^{\lambda }T^*_{X, J}$ par $F^{\lambda }_{_{J}}$. Les définitions précédentes sont compatibles avec les définitions classiques de la géométrie complexe. En effet dans le cas complexe intégrable les fibrés $F^{\lambda }_{_{J}}$ admettent une structure holomorphe canonique donné par le faisceau des sections holomorphes 
${\cal O}(F^{\lambda  }_{_{J}})$, qui est définie de façon naturelle à partir du faisceau  $\Omega^1 _{X,J}$. La connexion canonique de type $(0,1)$ sur le fibré $F^{\lambda }_{_{J}}$ induite par le faisceau ${\cal O}(F^{\lambda }_{_{J}})$ 
$$
\I_{_{{\cal O}( F^{\lambda }_{_{J}}  )}} \otimes_{_{{\cal O}_X } } \bar{\partial}_{_{J} }
$$
coïncide évidement avec la connexion $\bar{\partial}_{{F^{\lambda }_{_{J}} }}$ induite par la connexion 
$
\bar{\partial}_{_{J,1}}
$
et de plus on a toujours l'égalité évidente 
$$
{\cal O}(F^{\lambda }_{_{J}})=\Ker\,\bar{\partial}_{{F^{\lambda }_{_{J}} }}.
$$
Dans le cas d'une variété complexe intégrable on a donc le diagramme commutatif suivant.
\begin{diagram}[height=1cm,width=1cm]
\Omega^1 _{X,J}               &\longleftrightarrow &             \bar{\partial}_{_{J,1}}
\\
\dTo&                                                   &\dTo&  
\\
{\cal O}(F^{\lambda }_{_{J}})   &\longleftrightarrow &\bar{\partial}_{{F^{\lambda }_{_{J}} }}
\end{diagram}  
Dans le cas presque complexe non intégrable le faisceau ${\cal O}_X:=\Ker \bar{\partial}_{_{J}}$ est un faisceaux de fonctions constantes pour un choix générique de structure presque complexe $J$ non intégrable. Il suffit de prendre par example une structure fortement non intégrable. On rappelle qu'une structure presque complexe est dit fortement non-intégrable si le fibré tangent est engendré ponctuellement par les crochets des champs de vecteurs de type $(0,1)$. D'autre part dans cette situation il n'existe pas de repères locaux complexes $(\alpha  _k)_k\in {\cal E}(\Lambda ^{1,0}_{_J}T_X^*)^{\oplus n} (U)$ tels que 
$\bar{\partial}_{_{J,1}}\,\alpha _k\equiv 0$ sur l'ouvert $U$ pour tout $k=1,...,n$, car sinon ceci entraînerait que le fibré $\Lambda ^{1,0}_{_J}T_X^*$ est plat, ce qui n'est pas toujours le cas pour une variété presque complexe.  Un tel phénomène peut être aussi envisagé pour les fibrés 
$F^{\lambda }_{_{J}}$ et la connection canonique $\bar{\partial}_{{F^{\lambda }_{_{J}} }}$.
\\
Cependant la connexion $\bar{\partial}_{{F^{\lambda }_{_{J}} }}$ induit une structure holomorphe canonique sur toutes les restrictions 
$F^{\lambda }_{_{J}\,|_{\gamma (\Sigma)}}$ du fibré $F^{\lambda }_{_{J}}$ aux images des plongements $(j,J)$-holomorphes $\gamma :(\Sigma,j)\longrightarrow (X,J)$ d'une courbe holomorphe lisse $\Sigma\subset \C^m$. En effet la restriction
$$
\bar{\partial}_{ F^{\lambda }_{_{J}} |\gamma (\Sigma) }: {\cal E}\Big(F^{\lambda }_{_{J}\,|_{\gamma (\Sigma)}} \Big) \longrightarrow {\cal E}\Big(\Lambda ^{0,1}_{_J}T^*_{\gamma (\Sigma)}\otimes_{_{\C}} F^{\lambda }_{_{J}\,|_{\gamma (\Sigma)}}\Big)
$$
 de la connexion $\bar{\partial}_{{F^{\lambda }_{_{J}} }}$ est bien évidement intégrable étant donné que $\Lambda ^{0,2}_{_J}T^*_{\gamma (\Sigma)} =0$. La structure holomorphe canonique sur le fibré $F^{\lambda }_{_{J}\,|_{\gamma (\Sigma)}}$ est alors donné par la formule 
$$
{\cal O}\Big(F^{\lambda }_{_{J}\,|_{\gamma (\Sigma)}} \Big):=\Ker \Big(\bar{\partial}_{ F^{\lambda }_{_{J}} |\gamma (\Sigma) }\Big). 
$$
Soit $h$ une métrique hermitienne quelconque sur $F^{\lambda }_{_{J}}$. On définit alors la connexion de Chern $D^{h}_{{F^{\lambda  }_{_{J}} } }$ comme étant l'unique connexion hermitienne sur $F^{\lambda  }_{_{J}}$ telle que 
$$
(D^{h}_{{F^{\lambda }_{_{J}} } })^{0,1}=\bar{\partial}_{{F^{\lambda}_{_{J}} } }.
$$
\section{Expression locale des opérateurs $\partial_{_J },\,\bar{\partial}_{_J },\,\theta_{_J }$ et $\bar{\theta }_{_J }$.}
Soit $(\zeta _k)_{k} \in{\cal E}(T^{1,0}_{X,J})^{\oplus n}(U)$ un repère locale du fibré $T^{1,0}_{X,J}$ et $M^k,N^k,U^k,V^k\in M_n({\cal E}(U))$ 
les $n\times n$-matrices définies par les relations
\begin{align*}
[\bar{\zeta }_j,\bar{\zeta }_r]^{1,0}_{_{J} }=\sum_{k=1}^n\,N^k_{j,r}\,\zeta _k    &\qquad \qquad  
 [\bar{\zeta }_j,\bar{\zeta }_r]^{0,1} _{_{J} }=\sum_{k=1}^n\,M^k_{j,r}\,\bar{\zeta}_k \\
[ \zeta_j
 ,\zeta_r]^{1,0}_{_{J}}=\sum_{k=1}^n\,\overline{M}^k_{j,r}\,\zeta _k
&  \qquad\qquad  [\zeta _j,\zeta _r]^{0,1}_{_{J}}=\sum_{k=1}^n\,\overline{N}^k_{j,r}\,\bar{\zeta}_k \\
[\zeta_j,\bar{\zeta}_r]^{1,0}_{_{J}}=\sum_{k=1}^n\,U^k_{j,r}\,\zeta _k
 &\qquad\qquad  [\zeta_j,\bar{\zeta}_r]_{_{J}}^{0,1} =\sum_{k=1}^n\,V^k_{j,r}\,\bar{\zeta}_k
\end{align*}
On a les relations $M^k_{j,r}=-M^k_{r,j},\,N^k_{j,r}=-N^k_{r,j}$ et $V^k_{j,r}=-\overline{U}^k_{r,j}$. De plus on a l'expression locale
\begin{eqnarray*}
&\displaystyle{\tau_{_{J}}=\sum_{1\leq k<l\leq n}[\zeta _k,\zeta _l]^{0,1}_{_{J}}\otimes \zeta ^*_k\wedge\zeta ^*_l 
=\sum_{\substack{1\leq k<l\leq n
\\
\\
1\leq r\leq n }}\overline{N}^r_{k,l}\, \zeta ^*_k\wedge\zeta ^*_l \otimes \bar{\zeta}_r }&
\end{eqnarray*}
pour la forme de torsion de la structure presque complexe $J$. On rappelle que les éléments de l'espace vectoriel 
$\Lambda ^{p,q}_{_J}T_{X,x} ^*\otimes_{_{\C}} T_{X,J,x}$ s'identifient naturellement avec les éléments du type $u+\bar{u},\,u\in\Lambda^{p,q}_{_J}T_{X,x} ^*\otimes_{_{\C}} T^{1,0} _{X,J,x}$.
On introduit maintenant une notation très utile pour la suite. Soit $(\zeta _k)_{k} \in(T^{1,0}_{X,J,x})^{\oplus n}$ un repère. Alors 
$(\zeta _k+\bar{\zeta }_k)_{k} \in(T_{X,J,x})^{\oplus n}$ est un repère complexe de l'espace vectoriel $(T_{X,x},J_x)$. On notera  
$$
c\times_{_{J} }\zeta _k:=c\cdot \zeta _k+\bar{c}\cdot \bar{\zeta }_k
$$
l'opération de produit d'un scalaire $c\in \C$ avec le vecteur réel $\zeta _k+\bar{\zeta }_k\in T_{X,x}$. Si $\alpha \in \Lambda ^{p,q}_{_J}T_{X,x} ^*$ on notera 
$$
\alpha \otimes_{_{J}}\zeta _k:=\alpha \otimes \zeta _k+\bar{\alpha }\otimes \bar{\zeta }_k
$$ 
la $(p,q)$-forme à valeurs dans l'espace vectoriel $T_{X,J,x}$. Avec ces notations on aura par exemple l'expression locale suivante pour le tenseur de Nijenhuis
\begin{eqnarray*}
&\displaystyle{N_{_{J}}=\sum_{\substack{1\leq k<l\leq n
\\
\\
1\leq r\leq n }} N^r_{k,l}\,\bar{\zeta}^*_k\wedge \bar{\zeta}^*_l\otimes_{_{J}}\zeta _r=\sum_{1\leq k,l,r\leq n} N^r_{k,l}\,\bar{\zeta}^*_k\otimes \bar{\zeta}^*_l\otimes_{_{J}}\zeta _r}.& 
\end{eqnarray*}
Si $f$ est une fonction on a 
$\partial_{_J }f=\sum_{k=1}^n\,(\zeta_k\,. f)\,\zeta ^*_k\,,\; \bar{\partial}_{_J }f=\sum_{k=1}^n\,(\bar{\zeta}_k\,. f)\,\bar{\zeta}^*_k\,,\;\theta_{_J }f=0$ et $\bar{\theta }_{_J }f=0$.
De plus en utilisant les expressions intrinsèques des opérateurs $\partial_{_J },\,\bar{\partial}_{_J },\,\theta_{_J }$ et $\bar{\theta }_{_J }$ on a les expressions
\begin{align*}
\partial_{_J }\zeta _k^*=-\sum_{1\leq l<t\leq n}
\overline{M}^k_{l,t}\,\zeta _l^*\wedge\zeta _t^* &\qquad\qquad
\partial_{_J }\bar{\zeta}^*_k= \sum_{1\leq l,t\leq n} \overline{U}^k_{t,l}\,  \zeta _l^*\wedge\bar{\zeta}^*_t\\
\bar{\partial}_{_J }\zeta _k^*=-\sum_{1\leq l,t\leq n}U^k_{l,t} \,  \zeta _l^*\wedge\bar{\zeta}^*_t&\qquad\qquad
\bar{\partial}_{_J }\bar{\zeta}^*_k=-\sum_{1\leq l<t\leq n} M^k_{l,t}\,\bar{\zeta}_l^*\wedge\bar{\zeta}_t^*\\
\theta_{_J }\bar{\zeta}^*_k=  \sum_{1\leq l<t\leq n} \overline
{N}^k_{l,t}\,\zeta _l^*\wedge\zeta _t^*&\qquad\qquad
\bar{\theta }_{_J } \zeta ^*_k= \sum_{1\leq l<t\leq n} N^k_{l,t} \,\bar{\zeta}_l^*\wedge\bar{\zeta}_t^*     
\end{align*} 
Soit 
$$
u=\sum_{\substack{|K|=p
\\
\\
|L|=q} }\,u_{_{K,L} }\,\zeta ^*_{_{K} }\wedge  \bar{\zeta}^*_{_{L} } 
$$
une $(p,q)$-forme par rapport à la structure presque complexe $J$. Le fait que l'opérateur  $T:=\partial_{_J },\,\bar{\partial}_{_J },\,\theta_{_J }$ où $\bar{\theta }_{_J }$ vérifie la règle de Leibnitz implique l'égalité
$$
Tu=\sum_{\substack{|K|=p
\\
\\
|L|=q} }\Big(Tu_{_{K,L} }\wedge \zeta ^*_{_{K} }\wedge  \bar{\zeta}^*  _{_{L} }+
\sum_{j=1}^p(-1)^{j-1}u_{_{K,L} }T\zeta ^*_{k_j} \wedge\zeta ^*_{_{\hat{K}_j} }\wedge  \bar{\zeta}^*  _{_{L} }  
+
\sum_{j=1}^q(-1)^{p+j-1}u_{_{K,L} }T\bar{\zeta}^*_{l_j} \wedge\zeta ^*_{_{K} }\wedge  \bar{\zeta}^*  _{_{\hat{L}_j} }     \Big),
$$ 
où $\hat{K}_j:=(k_1,...,\hat{k}_j,...,k_p)$ et analoguement pour $\hat{L}_j$. On déduit alors les expressions locales
\begin{eqnarray}
&\displaystyle{\partial_{_J }u=\sum_{\substack{|K|=p
\\
\\
|L|=q} }\Big(\sum_{1\leq r\leq n}\,(\zeta _r\,.u_{_{K,L} })\,\zeta ^*_r\wedge\zeta ^*_{_{K} }\wedge  \bar{\zeta}^*  _{_{L} } +\sum _{\substack{1\leq j\leq p
\\
\\
1\leq r<t\leq n} }(-1)^{j}u_{_{K,L} } \cdot \overline{M}^{k_j}_{r,t}\,\zeta ^*_r\wedge\zeta ^*_t\wedge \zeta ^*_{_{\hat{K}_j} }\wedge  \bar{\zeta}^*  _{_{L} }            }&\nonumber
\\\nonumber
\\
&\displaystyle{-(-1)^p\sum _{\substack{1\leq j\leq q
\\
\\
1\leq r,t\leq n} }(-1)^{j}u_{_{K,L} } \cdot \overline{U}^{l_j}_{t,r}\,\zeta ^*_r\wedge\bar{\zeta}^*_t \wedge\zeta ^*_{_{K} }\wedge  \bar{\zeta}^*  _{_{\hat{L}_j} }    \Big)  }&\label{locdel} 
\\\nonumber
\\\nonumber
\\
&\displaystyle{\bar{\partial}_{_J }u=\sum_{\substack{|K|=p
\\
\\
|L|=q} }\Big(\sum_{1\leq r\leq n}\,(\bar{\zeta}_r\,.u_{_{K,L} })\,\bar{\zeta} ^*_r\wedge\zeta ^*_{_{K} }\wedge  \bar{\zeta}^*  _{_{L} } +\sum _{\substack{1\leq j\leq p
\\
\\
1\leq r,t\leq n} }(-1)^{j}u_{_{K,L} } \cdot U^{k_j}_{r,t}\,\zeta ^*_r\wedge\bar{\zeta}^*_t\wedge \zeta ^*_{_{\hat{K}_j} }\wedge  \bar{\zeta}^*  _{_{L} }            }&\nonumber
\\\nonumber
\\
&\displaystyle{+(-1)^p\sum _{\substack{1\leq j\leq q
\\
\\
1\leq r<t\leq n} }(-1)^{j}u_{_{K,L} } \cdot M^{l_j}_{r,t}\, \bar{\zeta}^*_r\wedge\bar{\zeta}^*_t \wedge\zeta ^*_{_{K} }\wedge  \bar{\zeta}^*  _{_{\hat{L}_j} }    \Big)  }&\label{locdelbar} 
\\\nonumber
\\
&\displaystyle{\theta_{_J } u=-(-1)^p\sum_{\substack{|K|=p
\\
\\
|L|=q} }\sum _{\substack{1\leq j\leq q
\\
\\
1\leq r<t\leq n} }(-1)^{j}u_{_{K,L} } \cdot \overline{N}^{l_j}_{r,t}\,\zeta ^*_r\wedge\zeta ^*_t\wedge \zeta ^*_{_{K} }\wedge  \bar{\zeta}^*  _{_{\hat{L}_j} }  }&\label{loctors} 
\end{eqnarray} 
\begin{eqnarray}
&\displaystyle{\bar{\theta }_{_J }u= -\sum_{\substack{|K|=p
\\
\\
|L|=q} }\sum _{\substack{1\leq j\leq p
\\
\\
1\leq r<t\leq n} }(-1)^{j}u_{_{K,L} } \cdot N^{k_j}_{r,t}  \, \bar{\zeta}^*_r\wedge\bar{\zeta}^*_t\wedge\zeta ^*_{_{\hat{K}_j} }\wedge  \bar{\zeta}^*  _{_{L} } }&\label{loctorsbar}
\end{eqnarray}  
\section{Relation entre la connexion de Chern du fibré tangent $T_{X,J}$ d'une variété presque complexe et la connection de Levi-Civita} 
Pour $p=1$ la définition \eqref{(0,1)conn-sur-(p,0)} de la connexion $\bar{\partial}_{_{J,1}}$ s'écrit sous la forme
$$
\bar{\partial}_{_{J,1}}\alpha (\eta)\cdot \xi =\eta\,.\,\alpha (\xi )-\alpha ([\eta,\xi ]^{1,0}) 
$$
pour tout $\alpha \in {\cal E}^{1,0}_{_{X,J}}(U),\,\eta \in{\cal E}(T^{0,1}_{X,J})(U) $ et $\xi \in{\cal E}(T^{1,0}_{X,J})(U)$. La connexion duale
$$
\bar{\partial}_{_{T^{1,0} _{X,J} } }: {\cal E}(T^{1,0} _{X,J}) \longrightarrow {\cal E}(\Lambda ^{0,1}_{_J}T^*_{X}\otimes_{_{\C}}T^{1,0} _{X,J})
$$
sur le fibré $T^{1,0} _{X,J}$, définie par la formule 
$$
(\bar{\partial}_{_{J,1}}\alpha)\cdot\xi =
\bar{\partial}_{_{J}}(\alpha \cdot \xi )-\alpha \cdot \bar{\partial}_{_{T^{1,0} _{X,J} } }\xi 
$$
vérifie alors l'identité 
$$
\bar{\partial}_{_{T^{1,0} _{X,J} } }\xi (\eta)=[\eta,\xi ]^{1,0}.
$$
Soit $(\zeta _k)_{k} \in{\cal E}(T^{1,0}_{X,J})^{\oplus n}(U)$ un repère local du fibré $T^{1,0}_{X,J}$ et $A''_{_{J} }=\sum_r(A''_{_{J} })^r\,\bar{\zeta}^*_r$ la forme de connexion de $\bar{\partial}_{_{T^{1,0} _{X,J} } }$ par rapport au repère en question. On a alors
$$
\bar{\partial}_{_{T^{1,0} _{X,J} }}\zeta _j(\bar{\zeta}_r)=-[\zeta _j,\bar{\zeta}_r]^{1,0}=
\sum_k(A''_{_{J} })_{k,j}(\bar{\zeta}_r)\zeta _k=-\sum_k\,U^k_{j,r}\,\zeta _k.   
$$
On déduit alors la formule $(A''_{_{J} })^r_{k,j}=-U^k_{j,r}$.
En utilisant l'isomorphisme $\C$-linéaire canonique du fibré $T^{1,0} _{X,J}$ avec le fibré tangent $T_{X,J}$ on déduit la  connexion de type $(0,1)$ canonique 
$$
\bar{\partial}_{_{T_{X,J}}}: {\cal E}(T_{X,J}) \longrightarrow {\cal E}(\Lambda ^{0,1}_{_J}T^*_{X}\otimes_{_{\C}}T_{X,J})
$$
du fibré tangent $T_{X,J}$. De façon explicite on a pour tout champ de vecteurs réels $\xi ,\,\eta\in{\cal E}(T_X) (U)$ l'expression suivante :
\begin{eqnarray*}
&\displaystyle{
\bar{\partial}_{_{T_{X,J}}}\xi (\eta)=\bar{\partial}_{_{T_{X,J}}}\xi (\eta^{0,1})=
[\eta^{0,1} \xi ^{1,0} ]^{1,0}+[\eta^{1,0} ,\xi^{0,1}]^{0,1}}&
\\
\\
&\displaystyle{=\frac{1}{4}\Big([\eta,\xi ]+[J\eta,J\xi ]+J[J\eta,\xi ]-J[\eta,J\xi ]\Big)}.& 
\end{eqnarray*}  
Soit $\omega \in {\cal E}(\Lambda ^{1,1}_{_J}T_X^*)(X)$ une métrique hermitienne sur $T_{X,J}$. On désignera par
$$
D^{\omega} _{_{J} }: {\cal E}(T_{X,J}) \longrightarrow {\cal E}(T^*_X\otimes_{_{\R}}T_{X,J})$$
la connexion de Chern du fibré  hermitien $(T_{X,J},\omega )$, autrement dit l'unique connexion $\omega$-hermitienne telle que 
$$
(D^{\omega} _{_{J}})^{0,1}=\bar{\partial}_{_{T_{X,J}}}.
$$ 
Considérons maintenant la métrique riemannienne $J$-invariante associée $g:=\omega (\cdot,J\cdot)\in{\cal E}(S^2_{_{\R} }T^*_X)(X)$. On désigne par
$$
\nabla^g:{\cal E}(T_X) \longrightarrow {\cal E}(T^*_X\otimes_{_{\R}}T_X)
$$
la connexion de Levi-Civita relative à la métrique riemannienne $g$.
Dans la suite on aura besoin de considérer la décomposition
$$
\Lambda ^k_{_{\R} } T_X^*\otimes_{_{\R}}T_{X,J} \simeq_{_{\C} }  \Lambda ^k_{_{\C} }(T_X\otimes_{_{\R}}\C)^*\otimes_{_{\C} }T_{X,J}\simeq_{_{\C} } \bigoplus_{p+q=k}   \Lambda^{p,q}_{_J}T_X^*\otimes_{_{\C}}T_{X,J}. 
$$
Le théorème suivant relie la connexion de Levi-Civita avec une connexion fondamentale de la géométrie presque complexe. Une autre formule peut être trouve dans \cite{Gau}. 
\begin{theo}\label{teorconchern}
Soit $(X,J)$ une variété presque complexe, $\omega \in {\cal E}(\Lambda ^{1,1}_{_J}T_X^*)(X)$ une métrique hermitienne sur $T_{X,J}$ et $g:=\omega (\cdot,J\cdot)\in{\cal E}(S^2_{_{\R} }T^*_X)(X)$ la métrique riemannienne $J$-invariante associée à $\omega$. Il existe deux tenseurs réels
$$
\delta_{_{J} } \omega   \in {\cal E}((T_X^*)^{\otimes 2} \otimes_{_{\R}}T_X)(X)
\quad
\mbox{et} 
\quad
N^{\omega} _{_J}\in {\cal E}((T^{0,1}_{X,J})^{*,\otimes 2} \otimes_{_{\C}} T_{X,J})(X)
$$
tels que $d\omega =0$ si et seulement si $\delta_{_{J} } \omega =0$; $N_{_{J} }=0 $ si et seulement si $N^{\omega} _{_J}=0$. La connexion de Chern $D^{\omega} _{_{J} }$ du fibré hermitien $(T_{X,J},\omega )$ est relié à la connection de Levi-Civita $\nabla^g$ par la formule
\begin{eqnarray}\label{Formu-Chern-LCivi}  
D^{\omega} _{_{J},\,\xi  }\,\eta:=\nabla^g_{\xi }\,\eta +\delta_{_{J} } \omega(\xi ,\eta)-N^{\omega} _{_J}(\xi ,\eta)
\end{eqnarray}
pour tout champ de vecteurs réels $\xi ,\,\eta\in{\cal E}(T_X) (U),\,(U\subseteq X \mbox{ouvert arbitraire}) $. Le $2$-tenseur réel $\delta_{_{J} } \omega$ est défini par la formule
$$
2\delta_{_{J} } \omega(\xi ,\eta):=\Big[\gamma^{2,0}_{_{\omega ,J}}+
\gamma^{0,2}_{_{\omega ,J}}\Big](\xi ,\eta)+J\gamma^{1,1} _{_{\omega ,J}}(\xi ,J\eta)
$$
où 
$$
\gamma^{2,0}_{_{\omega ,J}}\in {\cal E}(\Lambda^{2,0}_{_J}T_X^*\otimes_{_{\C}}T_{X,J} )(X),\;
\gamma^{1,1} _{_{\omega ,J}}\in {\cal E}(\Lambda ^{1,1}_{_J}T_X^*\otimes_{_{\C}}T_{X,J}  )(X)\;\mbox{et} \; 
\gamma^{0,2}_{_{\omega ,J}}\in {\cal E}(\Lambda^{0,2}_{_J}T_X^*\otimes_{_{\C}}T_{X,J} )(X)
$$ 
sont les composantes, $($par rapport à la structure presque complexe $J$$)$ de la $2$-forme réelle 
\\
$\gamma_{_{\omega}}\in {\cal E}(\Lambda ^2T_X^*\otimes_{_{\R}}T_X)(X)$ définie par la formule 
$$
\omega (\gamma_{_{\omega}}(\xi ,\eta),\mu )=d\omega (\xi ,\eta,\mu )
$$
pour tout champ de vecteurs réels
$\xi ,\eta,\mu\in {\cal E}(T_X)(X)$. Enfin le $(0,2)$-tenseur réel $N^{\omega} _{_J}$ est définie par la formule $N^{\omega} _{_J}:=\tau^{\omega} _{_J}+\bar{\tau}^{\omega} _{_J}$ où
$$
\tau^{\omega} _{_J}\in {\cal E}((T^{0,1}_{X,J})^{*,\otimes 2} \otimes_{_{\C}} T^{1,0}_{X,J})(X)
$$ 
est
le $(0,2)$-tenseur défini par la formule
$$
\omega (\tau^{\omega} _{_J}(\xi ,\eta),\mu )=\omega (\xi ,[\eta,\mu]^{1,0})
$$
pour tout $(0,1)$-champ de vecteurs $\xi ,\,\eta,\,\mu \in{\cal E}(T^{0,1}_{X,J})(X)$. La forme de torsion ${\cal T} _{D^{\omega} _{_{J} }}$ de la connexion de Chern $D^{\omega} _{_{J} }$ vérifie l'identité
\begin{eqnarray}\label{For-Tor-Ch-LC}  
{\cal T} _{D^{\omega} _{_{J} }}(\xi ,\eta)=
-\Big[\gamma^{2,0}_{_{\omega ,J}}+\gamma^{0,2}_{_{\omega ,J}}\Big](\xi ,\eta)-N^{\omega} _{_J}(\xi ,\eta)+N^{\omega} _{_J}(\eta,\xi ).
\end{eqnarray}
Si $N_{_{J} }=0$ alors $\gamma^{0,2}_{_{\omega ,J}}=0$. 
\end{theo}  
$Preuve$
\\
{\bf Expression  de la connexion de Chern $D^{\omega} _{_{J}}$ du fibré hermitien $(T_{X,J},\omega )$.}
\\
Soit $h_{\omega}$ la forme hermitienne sur le fibré $T_{X,J}$ associée à $\omega$. On rappelle qu'elle est définie par la formule $h_{\omega}(\xi ,\eta):=\omega (\xi ,J\eta)-i\omega (\xi ,\eta)$. La connexion de Chern $D^{\omega} _{_{J}}$ est définie par les formules
\begin{eqnarray}\label{Def-Chern}
&\displaystyle{
D^{\omega} _{_{J},\,\xi  }\,\eta=D^{\omega} _{_{J},\,\xi^{1,0}}\,\eta+\bar{\partial}_{_{T_{X,J}}}\eta (\xi^{0,1}),}&\nonumber
\\\nonumber
\\ 
&\displaystyle{
h_{\omega}(D^{\omega} _{_{J},\,\xi^{1,0}}\,\eta\,,\,\mu)=\xi ^{1,0}.\,h_{\omega} (\eta,\mu)-h_{\omega}(\eta,\bar{\partial}_{_{T_{X,J}}}\eta (\xi^{0,1}))}&  
\end{eqnarray}
pour tout champ de vecteurs réels $\xi ,\,\eta,\,\mu \in{\cal E}(T_X)(U)$. L'identité
$h_{\omega}(\xi ,\eta)=h_{\omega}(\xi^{1,0}  ,\eta^{0,1} )=-2i\omega(\xi^{1,0}  ,\eta^{0,1} )$ et la définition de la connexion canonique $\bar{\partial}_{_{T_{X,J}}}$ montrent que la formule \ref{Def-Chern} est équivalente à la formule
$$
\omega(D^{\omega} _{_{J},\,\xi^{1,0}}\,\eta\,,\,\mu^{0,1} )=\xi ^{1,0}.\,\omega (\eta^{1,0},\mu ^{0,1})
-\omega (\eta^{1,0},[\xi ^{1,0}, \mu ^{0,1}]^{0,1}) 
$$
On obtient en conclusion que la connexion de Chern peut étre définie par la formule
\begin{eqnarray}\label{Def-Chern-Expl}
\omega(D^{\omega} _{_{J},\,\xi}\,\eta\,,\,\mu^{0,1} )=\xi ^{1,0}.\,\omega (\eta^{1,0},\mu ^{0,1})
-\omega (\eta^{1,0},[\xi ^{1,0}, \mu ^{0,1}]^{0,1})+\omega ([\xi ^{0,1},\eta^{1,0}]^{1,0},\mu^{0,1})
\end{eqnarray}
pour tout champ de vecteurs réels $\xi ,\,\eta,\,\mu \in{\cal E}(T_X)(U)$.
\\
\\
{\bf{Expression de la connexion de Levi-Civita $\nabla^{g}$.}} 
\\
La  connexion de Levi-Civita
$
\nabla^g:{\cal E}(T_X) \longrightarrow {\cal E}(T^*_X\otimes_{_{\R}}T_X)
$
est définie par la formule classique
\begin{eqnarray*} 
&\displaystyle{
2g(\nabla^g_{\xi }\,\eta,\mu )=\xi \,.g(\eta,\mu )-\mu \,.g(\xi ,\eta)+\eta\,.g(\mu ,\xi)}&
\\
\\
&\displaystyle{-g(\xi ,[\eta,\mu ])+g(\mu ,[\xi ,\eta])+g(\eta,[\mu ,\xi ])}& 
\end{eqnarray*}  
pour tout champ de vecteurs réels $\xi, \,\eta,\,\mu \in {\cal E}(T_X)(U)$. Bien évidemment la définition précédente est équivalente à la formule
\begin{eqnarray}\label{Levi-Civita} 
&\displaystyle{
2\omega (\nabla^g_{\xi }\,\eta,-i\mu^{0,1})=\xi \,.\,\omega (\eta^{1,0} ,-i\mu^{0,1})
-\mu^{0,1} .\,\omega (\xi ,J\eta)+\eta\,.\,\omega (\mu^{0,1},i\xi^{1,0} )}&\nonumber
\\\nonumber
\\
&\displaystyle{-\omega(\xi ,J[\eta,\mu^{0,1}])+\omega(\mu^{0,1},i[\xi ,\eta]^{1,0})
+\omega(\eta,J[\mu^{0,1},\xi]).}& 
\end{eqnarray}
\\
{\bf{Expression des $2$-tenseurs 
$\gamma^{2,0} _{_{\omega ,J}},\,\gamma^{1,1} _{_{\omega ,J}}(\cdot , J\cdot)$ et $\gamma^{0,2}_{_{\omega ,J}}$.} } 
\\
On rappelle que les éléments de l'espace vectoriel 
$\Lambda ^{p,q}_{_J}T_{X,x} ^*\otimes_{_{\C}} T_{X,J,x}$ s'identifient naturellement avec les éléments du type $u+\bar{u},\,u\in\Lambda ^{p,q}_{_J}T_{X,x} ^*\otimes_{_{\C}} T^{1,0} _{X,J,x}$, (voir la section \ref{Con-presq-compl}). On a donc les identitées
$$
\gamma^{2,0} _{_{\omega ,J}}=\hat{\gamma}^{2,0} _{_{\omega ,J}}+\overline{\hat{\gamma}^{2,0} _{_{\omega ,J}}} ,
\quad
\gamma^{1,1} _{_{\omega ,J}}=\hat{\gamma}^{1,1} _{_{\omega ,J}} +\overline{\hat{\gamma}^{1,1} _{_{\omega ,J}}},
\quad 
\gamma^{0,2}_{_{\omega ,J}}=\hat{\gamma}^{0,2}_{_{\omega ,J}}+\overline{\hat{\gamma}^{0,2}_{_{\omega ,J}}}, 
$$
sur $T_X$, avec 
$$
\hat{\gamma}^{2,0}_{_{\omega ,J}}\in {\cal E}(\Lambda^{2,0}_{_J}T_X^*\otimes_{_{\C}}T^{1,0} _{X,J} )(X),\;
\hat{\gamma}^{1,1} _{_{\omega ,J}}\in {\cal E}(\Lambda ^{1,1}_{_J}T_X^*\otimes_{_{\C}}T^{1,0}_{X,J}  )(X)\;
\mbox{et} 
\; 
\hat{\gamma}^{0,2}_{_{\omega ,J}}\in {\cal E}(\Lambda^{0,2}_{_J}T_X^*\otimes_{_{\C}}T^{1,0}_{X,J} )(X).
$$ 
La décomposition 
$d\omega =\partial_{_J }\omega +\bar{\partial}_{_J }\omega -\theta_{_J }\omega -
\bar{\theta }_{_J }\omega $ implique alors les identitées
\begin{eqnarray*}
&\displaystyle{
\omega (\hat{\gamma}^{2,0}_{_{\omega ,J}}(\xi,\eta ),\mu^{0,1})=
\omega (\hat{\gamma}^{2,0}_{_{\omega ,J}}(\xi^{1,0} ,\eta^{1,0} ),\mu^{0,1})=
\partial_{_J }\omega(\xi^{1,0} ,\eta^{1,0} ,\mu^{0,1}),}&
\\
\\
\\
&\displaystyle{
\omega (\hat{\gamma}^{1,1}_{_{\omega ,J}}(\xi,J\eta ),\mu^{0,1})=
\omega (\hat{\gamma}^{1,1}_{_{\omega ,J}}(\xi^{1,0} ,-i\eta^{0,1}),\mu^{0,1})+
\omega (\hat{\gamma}^{1,1}_{_{\omega ,J}}(\xi^{0,1} ,i\eta^{1,0}),\mu^{0,1})=
}&\nonumber
\\\nonumber
\\
&\displaystyle{
=\bar{\partial}_{_J }\omega(\xi^{1,0} ,-i\eta^{0,1},\mu^{0,1})+
\bar{\partial}_{_J }\omega(\xi^{0,1} ,i\eta^{1,0},\mu^{0,1}),
}&
\\
\\
\\
&\displaystyle{
\omega (\hat{\gamma}^{0,2}_{_{\omega ,J}}(\xi^{0,1} ,\eta^{0,1}),\mu^{0,1})=
-\bar{\theta }_{_J }\omega(\xi^{0,1} ,\eta^{0,1},\mu^{0,1})
}&
\end{eqnarray*}
En explicitant les formes $\partial_{_J }\omega,\,\bar{\partial}_{_J }\omega $ et 
$\bar{\theta }_{_J }\omega$ dans les identitées précédentes on obtient les expressions suivantes
\begin{eqnarray} \label{Form-2,0} 
&\displaystyle{
\omega (\hat{\gamma}^{2,0}_{_{\omega ,J}}(\xi,\eta ),\mu^{0,1})
=\xi ^{1,0}.\,\omega (\eta^{1,0},\mu ^{0,1})-\eta^{1,0}.\,\omega (\xi^{1,0},\mu ^{0,1}) 
}&\nonumber
\\\nonumber
\\
&\displaystyle{
-\omega ([\xi^{1,0},\eta^{1,0}]^{1,0},\mu ^{0,1})  +
\omega ([\xi^{1,0},\mu ^{0,1}]^{0,1},\eta^{1,0} )-
\omega ([\eta^{1,0},\mu ^{0,1} ]^{0,1},\xi^{1,0} ),                 
}&
\end{eqnarray} 
\\
\begin{eqnarray}\label{Form-1,1}  
&\displaystyle{
\omega (\hat{\gamma}^{1,1}_{_{\omega ,J}}(\xi,J\eta ),\mu^{0,1})=
i\eta^{0,1}.\,\omega (\xi ^{1,0},\mu ^{0,1})-
i\mu ^{0,1}.\,\omega (\xi ^{1,0},\eta^{0,1})
}&
\nonumber
\\\nonumber
\\
&\displaystyle{
+i\xi ^{0,1}.\,\omega (\eta^{1,0},\mu ^{0,1}) -
i\mu ^{0,1}.\,\omega (\eta^{1,0},\xi ^{0,1})   
}&\nonumber
\\\nonumber
\\
&\displaystyle{
+i\omega ([\xi^{1,0},\eta^{0,1}]^{1,0} ,\mu^{0,1})-
i\omega ([\xi^{1,0},\mu ^{0,1}]^{1,0} ,\eta^{0,1})+
i\omega ([\eta^{0,1},\mu ^{0,1}]^{0,1} ,\xi^{1,0})
}&\nonumber
\\\nonumber
\\
&\displaystyle{
+i\omega ([\eta^{1,0},\xi^{0,1}]^{1,0} ,\mu^{0,1})-
i\omega ([\eta^{1,0},\mu^{0,1}]^{1,0} ,\xi^{0,1})+
i\omega ([\xi ^{0,1},\mu ^{0,1}]^{0,1} ,\eta^{1,0})
}&
\end{eqnarray}
et en fin
\begin{eqnarray}\label{Form-0,2} 
&\displaystyle{
\omega (\hat{\gamma}^{0,2}_{_{\omega ,J}}(\xi^{0,1} ,\eta^{0,1}),\mu^{0,1})=
-
\omega([\xi^{0,1},\eta^{0,1}]^{1,0},\mu^{0,1})
}&\nonumber
\\\nonumber
\\
&\displaystyle{
+
\omega([\xi^{0,1},\mu ^{0,1}]^{1,0},\eta^{0,1})-
\omega([\eta^{0,1},\mu^{0,1}]^{1,0},\xi^{0,1}).}&
\end{eqnarray}
En remplaçant $-i\mu ^{0,1}$ à la place de $\mu ^{0,1}$ dans les identitées \ref{Form-2,0} et \ref{Form-0,2}, en sommant les identitées obtenues avec l'identité \ref{Form-1,1} et en tenant compte de la formule \ref{Def-Chern-Expl} on obtient l'identité voulue
\ref{Formu-Chern-LCivi}. Le fait que le $2$-tenseur $\gamma^{1,1} _{_{\omega ,J}}(\cdot , J\cdot)$ est symétrique implique l'identité \ref{For-Tor-Ch-LC} sur la forme de torsion de la connexion de Chern.\hfill$\Box$
\section{La courbure de Chern des puissances de Schur du fibré des $(1,0)$-formes}
On a la définition suivante.
\begin{defi}
Le tenseur de courbure de Chern 
$$
{\cal C}_h(F^{\lambda }_{_{J}})
\in {\cal E}(\Lambda ^{1,1}_{_{J}}T^*_X
\otimes_{_{\C} } F^{\lambda,\,* }_{_{J}}\otimes_{_{\C} }F^{\lambda  }_{_{J}})(X)  
$$ 
du fibré vectoriel hermitien $(F^{\lambda }_{_{J}},h)\longrightarrow (X,J)$ est la $(1,1)$-forme donnée par la formule 
$$
{\cal C}_h(F^{\lambda  }_{_{J}}):=\Theta (D^{h}_{{F^{\lambda }_{_{J}} } })^{1,1}.
$$
La courbure de Chern 
$$
{\cal C}^h_{F^{\lambda }_{_{J}}}\in
 {\cal E}(\Herm(T_{X,J} \otimes_{_{\C}}F^{\lambda  }_{_{J}}))(X)   
$$
est la forme hermitienne sur le fibré vectoriel complexe $T_{X,J} \otimes_{_{\C}}F^{\lambda  }_{_{J}}$ définie par la formule 
$$
{\cal C}^h_{F^{\lambda }_{_{J}}}(\xi \otimes \sigma ,\eta\otimes \tau)
:=h({\cal C}_h(F^{\lambda }_{_{J}})(\xi ^{1,0}_{_{J}},\eta^{0,1}_{_{J}} )\cdot \sigma ,\tau)
$$
pour tout champ de vecteurs réels $\xi ,\eta\in {\cal E}(T_X)(U)$ et sections $\sigma,\tau \in {\cal E}(F^{\lambda  }_{_{J}})(U)$ sur un ouvert $U$ quelconque.
\end{defi} 
La courbure de Chern
${\cal C}^h_{F^{\lambda  }_{_{J}}}$  est une forme hermitienne sur le fibré vectoriel complexe $T_{X,J} \otimes_{_{\C}}F^{\lambda  }_{_{J}}$ grâce à la relation $\eqref{hermitcurvat}$ (remarquée dans la section \ref{basic-herm-con}). Soit
$$
{\cal C}_h(F^{\lambda }_{_{J}})=\sum_{\substack{1\leq l,m \leq r_{\lambda } 
\\
\\
1\leq j,k\leq n }} C^{j,k}_{l,m}\,\zeta ^*_j\wedge \bar{\zeta } ^*_k\otimes  e^*_m\otimes e_l
$$
l'expression locale du tenseur de courbure de Chern, (ici $r_{\lambda }:=rg_{_{\C} }F^{\lambda }_{_{J}}$). Si le repère local $(e_l)_{l}\in {\cal E}(F^{\lambda }_{_{J}})^{\oplus r_{\lambda } } (U)$ est $h(x_0)$-orthonormé en un point $x_0$ alors l'expression locale de la courbure de Chern s'écrit sous la forme
$$
{\cal C}^h_{F^{\lambda }_{_{J}}}(x_0) =\sum_{\substack{1\leq l,m \leq r_{\lambda } 
\\
\\
1\leq j,k\leq n }} C^{j,k}_{l,m}(x_0)\,\zeta ^*_j\otimes  e^*_m\otimes \bar{\zeta } ^*_k\otimes \bar{e}^*_l
$$
où les coefficients vérifient la relation $\overline{C^{j,k}_{l,m}(x_0)}=C^{k,j}_{m,l}(x_0)$ vue dans la section \ref{basic-herm-con}. Remarquons que $(\zeta^*_{k|_{T_X}})_k\in {\cal E}(T^*_{X,J} )^{\oplus n} (U)$ est le repère dual du repère 
$(\zeta _k+\bar{\zeta }_k )_k\in {\cal E}(T_{X,J} )^{\oplus n} (U)$ par rapport à la structure $J$. Bien évidemment il est équivalent de donner soit le tenseur de courbure soit la courbure de Chern. On aura besoin de la définition suivante.
\begin{defi}
Une section $\sigma \in {\cal E}(F^{\lambda }_{_{J}})(U)$ est dite presque-holomorphe au point $x\in U$ si on a
$\bar{\partial}\,\sigma (x)=0$. 
Un repère local $(\sigma_k)_k \subset {\cal E}(F^{\lambda }_{_{J}})(U)$ est dit 
presque-holomorphe spécial au point $x\in U$ si 
$\bar{\partial}\,\sigma_k (x)=0$ et 
$(D^h)^{1,0} \,\bar{\partial}\,\sigma_k (x)=0$ pour tout $k$.
\end{defi}  
La définition de repère local presque-holomorphe spécial en un point est indépendante de la métrique hermitienne. En effet si $A''_{\sigma } $ est la matrice de la connexion de type $(0,1)$ canonique du fibré vectoriel $F^{\lambda }_{_{J}}$ relative au repère $(\sigma_k)_k \subset {\cal E}(F^{\lambda }_{_{J}})(U)$, la condition que le repère local $(\sigma_k)_k$ soit presque-holomorphe spécial au point $x$ s'exprimé par les égalités $A''_{\sigma } (x)=0$ et $\partial_{_{J} } A''_{\sigma } (x)=0$.
Le lemme élémentaire suivant donne une première idée de l'utilité de la notion de courbure de Chern.
\begin{lem}\label{lemchern}
Soient $\sigma ,\,\tau\in {\cal E}(F^{\lambda }_{_{J}})(U)$ deux sections presque-holomorphes en un point $x\in U$ du fibré hermitien $(F^{\lambda }_{_{J}},h)\longrightarrow (X,J)$ et
$\xi ,\eta\in {\cal E}(T_X)(U)$ deux champs de vecteurs réels. Alors au point $x$ on a l'identité
\begin{eqnarray}\label{chernscal}
&\displaystyle{
{\cal C}^h_{F^{\lambda }_{_{J}}}(\xi \otimes \sigma ,\eta\otimes \tau)_{|x} 
=\bar{\partial}_{_{J}}\partial_{_{J}}\,h(\sigma ,\tau) (\xi^{1,0},\eta^{0,1})_{|x} +
h(\xi ^{1,0}_{_{D} }.\,\sigma ,\eta^{1,0}_{_{D} }.\,\tau)_{|x} }&\nonumber
\\\nonumber
\\
&\displaystyle{+h(\xi ^{1,0}_{_{D} }.\,\eta^{0,1}_{_{D} }.\,\sigma ,\tau)_{|x} 
+h(\sigma ,\eta^{1,0}_{_{D} }.\,\xi^{0,1}_{_{D} }.\,\tau )_{|x} }.&
\end{eqnarray}
Soit $(\sigma_k)_k \subset {\cal E}(F^{\lambda }_{_{J}})(U)$ un repère local
presque-holomorphe spécial au point $x\in U$. Alors au point $x$ on a l'identité
\begin{eqnarray}\label{chernscalspec}
{\cal C}^h_{F^{\lambda }_{_{J}}}(\xi \otimes \sigma_k ,\eta\otimes \sigma _l)_{|x} 
=\bar{\partial}_{_{J}}\partial_{_{J}}\,h(\sigma_k ,\sigma _l) (\xi^{1,0},\eta^{0,1})_{|x} +
h(\xi ^{1,0}_{_{D} }.\,\sigma_k ,\eta^{1,0}_{_{D} }.\,\sigma _l)_{|x}.
\end{eqnarray}
En particulier
\begin{eqnarray}\label{chernscalpsh}
i\partial_{_{J}} \bar{\partial}_{_{J}}\,|\sigma_k |^2_h\,(\xi ,J\xi )_{|x}
=-2\,{\cal C}^h_{F^{\lambda }_{_{J}}}\,(\xi \otimes \sigma_k ,\xi \otimes \sigma_k)_{|x}+2\,|\xi ^{1,0}_{_{D} }.\,\sigma_k|^2_{h\,|x}.
\end{eqnarray}
\end{lem} 
Dans le cas d'une variété complexe $(X,J)$ et d'un fibré vectoriel holomorphe hermitien $(F,h)\longrightarrow (X,J)$ on a pour toutes sections holomorphes $\sigma ,\,\tau \in {\cal O}(F)(U)$ l'identité
$$
{\cal C}^h_{F}(\xi \otimes \sigma ,\eta\otimes \tau)
=\bar{\partial}_{_{J}}\partial_{_{J}}\,h(\sigma ,\tau) (\xi^{1,0},\eta^{0,1})+
h(\xi ^{1,0}_{_{D} }.\,\sigma ,\eta^{1,0}_{_{D} }.\,\tau)
$$
sur l'ouvert $U$. On déduit en particulier la formule remarquable suivante
$$
i\partial_{_{J}} \bar{\partial}_{_{J}}\,|\sigma |^2_h\,(\xi ,J\xi )
=-2\,{\cal C}^h_{F}\,(\xi \otimes \sigma ,\xi \otimes \sigma)+2\,|\xi ^{1,0}_{_{D} }.\,\sigma|^2_h
$$
qui montre que pour tout section holomorphe $\sigma  \in {\cal O}(F)(U)$ la fonction $|\sigma |^2_h$ est plurisousharmonique sur l'ouvert $U$ si la courbure du fibré $F$ est négative au sens de Griffiths, autrement dit si 
${\cal C}^h_{F}(\xi \otimes \sigma ,\xi \otimes \sigma)\leq 0$ pour tout $\xi \in T_{X,x}$ et $\sigma \in F_x$, (voir \cite{Gri} et \cite{Dem-1}, chapitre VII pour des applications fondamentales de la notion de courbure au sens de Griffiths). On déduit en particulier que si la variété complexe $X$ est compacte, connexe et $\sigma  \in {\cal O}(F)(X)$ est une section globale d'un fibré vectoriel holomorphe $F$ admettant une métrique hermitienne à courbure négative au sens de Griffiths alors le section $\sigma $ est identiquement nulle sur $X$ si elle s'annule en un point. On remarque que la notion de positivité (négativité) au sens de Griffiths pour un fibré $(F^{\lambda }_{_{J}},h)$ ne signifie rien d'autre que pour tout vecteur réel $\xi \in T_{X,J}$ l'endomorphisme $h$-hermitien $i{\cal C}_h(F^{\lambda }_{_{J}})(\xi ,J\xi )$ est positif (négatif). Si la courbure du fibré $(F^{\lambda }_{_{J}},h)$ est strictement négative au sens de Griffiths en un point $x$ alors on déduit d'après la formule $\eqref{chernscalpsh}$ que les fonctions $|\sigma_k |^2_h$ sont strictement $J$-plurisousharmoniques au voisinage du point $x$, (voir \cite{Pal} pour la notion de fonction strictement $J$-plurisousharmoniques et pour plus de détails). 
\\
\\
$Preuve\;du\;lemme\;\ref{lemchern} $. On a l'égalité
$$
\partial_{_{J}} \bar{\partial}_{_{J}}\,h(\sigma ,\tau)=\{D^{1,0}_{_{F} } \,\bar{\partial}_{_{F} }\sigma ,\tau \}_h-
\{\bar{\partial}_{_{F} }\sigma,\bar{\partial}_{_{F} }\tau\}_h+ \{D^{1,0}_{_{F} }\sigma ,D^{1,0}_{_{F} }\tau\}_h +\{\sigma ,\bar{\partial}_{_{F} }D^{1,0}_{_{F} }\tau\}_h.
$$
Le fait que $\deg\,\sigma =\deg\,\tau=0$ et l'identité 
$\{{\cal C}_h(F)\cdot \sigma ,\tau\}_h+\{\sigma ,{\cal C}_h(F)\cdot\tau\}_h=0$ impliquent
\begin{eqnarray}\label{orror} 
&\displaystyle{
\partial_{_{J}} \bar{\partial}_{_{J}}\,h(\sigma ,\tau)=-\{{\cal C}_h(F)\cdot \sigma ,\tau\}_h-
\{\sigma ,D^{1,0}_{_{F} }\bar{\partial}_{_{F} }\tau\}_h+ \{D^{1,0}_{_{F} } \,\bar{\partial}_{_{F} }\sigma ,\tau \}_h}&\nonumber
\\\nonumber
\\
&\displaystyle{
-
\{\bar{\partial}_{_{F} }\sigma,\bar{\partial}_{_{F} }\tau\}_h+ \{D^{1,0}_{_{F} }\sigma ,D^{1,0}_{_{F} }\tau\}_h}.& 
\end{eqnarray}
En explicitant l'égalité précédente par rapport au champs de vecteurs réels $\xi$ et $\eta$ on obtient l'identité
\begin{eqnarray*}
&\displaystyle{
{\cal C}^h_{_{F} }(\xi \otimes \sigma ,\eta\otimes \tau) 
=\bar{\partial}_{_{J}}\partial_{_{J}}\,h(\sigma ,\tau) (\xi^{1,0},\eta^{0,1}) 
+h(\xi ^{1,0}_{_{D} }.\,\eta^{0,1}_{_{D} }.\,\sigma ,\tau)
+h(\sigma ,\eta^{1,0}_{_{D} }.\,\xi^{0,1}_{_{D} }.\,\tau )
 }&
\\
\\
&\displaystyle{+
h([\eta^{0,1} ,\xi ^{1,0} ]^{0,1}_{_{D} }.\,\sigma ,\tau)+h(\sigma ,[\xi ^{0,1} ,\eta^{1,0} ]^{0,1}_{_{D} }.\,\tau)+
h(\xi ^{1,0}_{_{D} }.\,\sigma ,\eta^{1,0}_{_{D} }.\,\tau)+h(\eta ^{0,1}_{_{D} }.\,\sigma ,\xi ^{0,1}_{_{D} }.\,\tau) }&
\end{eqnarray*} 
qui permet de déduire la formule $\eqref{chernscal}$. Soit $(\sigma _k)_k$ le repère de l'énonce du lemme. On déduit d'après l'identité $\eqref{orror}$ l'égalité suivante au point $x$;
$$
\partial_{_{J}} \bar{\partial}_{_{J}}\,h(\sigma_k ,\sigma _l)_{|x}=-\{{\cal C}_h(F)\cdot \sigma_k ,\sigma _l\}_{h\,|x}+\{D^{1,0}_{_{F} }\sigma_k ,D^{1,0}_{_{F} }\sigma _l\}_{h\,|x}
$$
qui permet de conclure la preuve du lemme.
\hfill$\Box$ 
\\
Dans la sous-section suivante on montre l'existence de repères locaux $(\sigma_k)_k \subset {\cal E}(F^{\lambda }_{_{J}})(U)$ presque-holomorphes spéciaux en un point $x\in U$ tels que $D^{\omega }_{_{J} }\sigma _k(x)=0$ pour tout $k$. Dans ce cas on déduit d'après les formules $\eqref{chernscalspec}$ et $\eqref{chernscalpsh}$ les identités suivantes au point $x$;\begin{eqnarray*}
&\displaystyle{
{\cal C}^{\omega}_{_{{X,J} } }(\xi \otimes \sigma_k ,\eta\otimes \sigma _l)_{|x} 
=\bar{\partial}_{_{J}}\partial_{_{J}}\,h (\sigma_k ,\sigma _l) (\xi^{1,0},\eta^{0,1})_{|x} }&
\\
\\
&\displaystyle{
i\partial_{_{J}} \bar{\partial}_{_{J}}\,|\sigma_k |^2_{h} \,(\xi ,J\xi )_{|x}
=-2\,{\cal C}^h_{F^{\lambda }_{_{J}}} \,(\xi \otimes \sigma_k ,\xi \otimes \sigma_k)_{|x},}&\end{eqnarray*}
pour tout champs de vecteurs réels $\xi ,\eta\in{\cal E}(T_X)(U)$.

\subsection{Interprétation géométrique de la notion de courbure de Chern dans le cas presque complexe}\label{intercurvat}
Le lemme fondamental suivant est une version presque complexe d'un lemme classique de la géométrie hermitienne complexe (voir \cite{Dem-1}, chapitre V).
\begin{lem} \label{senscurvat}
Soit $(X,J)$ une variété presque complexe et $(F^{\lambda }_{_{J}},h)\longrightarrow (X,J)$ le fibré vectoriel hermitien d'une puissance de Schur du fibré des $(1,0)$-formes. Soient $(z_1,...,z_n)$ des coordonnées $\ci$ complexes  centrées en un point $x$ telles que $J(x)=J_0$, où $J_0$ désigne la structure presque complexe canonique relative à ces coordonnées. Il existe un repère local $(\sigma _k)_k\in {\cal E}(F^{\lambda }_{_{J}})^{\oplus r_{\lambda} } (U_x)$ presque-holomorphe spécial au point $x$ pour lequel les coefficients de la métrique hermitienne $h$ s'écrivent sous la forme
$$
h(\sigma _l,\sigma_m)=\delta _{l,m}+\sum_{1\leq j,k\leq n}\,H^{j,\bar{k}}_{l,m}\,z_j\bar{z}_k+O(|z|^3). $$
Quel que soit le choix du repère $(\sigma _k)_k\in {\cal E}(F^{\lambda }_{_{J}})^{\oplus r_{\lambda} } (U_x)$ presque-holomorphe spécial au point $x$ pour lequel les coefficients de la métrique hermitienne $h$ s'écrivent sous la forme précédente on a les expressions suivantes pour le tenseur de courbure et la courbure de Chern au point $x$:
\begin{eqnarray}
&\displaystyle{
{\cal C}_h(F^{\lambda  }_{_{J}})_{|x} =-\sum_{ \substack{1\leq l,m\leq r_{\lambda}
\\
\\
1\leq j,k\leq n }  } H^{j,\bar{k}}_{l,m} \,dz_j\wedge d\bar{z}_k\otimes \sigma  ^*_m\otimes\sigma _l}&\label{curvat1} 
\\\nonumber
\\
&\displaystyle{
{\cal C}^h_{F^{\lambda }_{_{J}}} (\xi   \otimes \sigma _l ,\eta\otimes \sigma _m)_{|x} 
=\bar{\partial}_{_{J}}\partial_{_{J}}\,h (\sigma_l ,\sigma _m) (\xi ^{1,0},\eta^{0,1})_{|x} },&
\label{curvat2}
\end{eqnarray}
pour tout champ de vecteurs réels $\xi  ,\eta \in{\cal E}(T_X)(U_x)$ et tout indice $l,\,m$.
\end{lem} 
Le lemme nous montre que la courbure de Chern au point $x$ mesure l'obstruction à l'existence de repères locaux presque-holomorphes spéciaux et orthonormaux à l'ordre deux en $x$. 
\\
\\
$Preuve$.
Soit $e\equiv (e_k)_k\in {\cal E}(F^{\lambda }_{_{J}})^{\oplus r_{\lambda} }(U_x)$ un repère local $h(x)$-orthonormé au point $x$. On peut supposer que la forme de la connexion de Chern 
$D^h_{F^{\lambda }}$ relative à ce repère vérifie la condition $A_{e}(x)=0$. En effet en effectuant un changement de repère $e'=e\cdot g_0$ avec $g_0=\I+O(|z|),\, dg_0(x)=-A_{e}(x)$ on a que la forme de connexion $A_{e'}=g_0^{-1}(dg_0+A_{e}\cdot g_0)$ relative au repère $e'$ vérifie la propriété voulue. Soient
\begin{eqnarray*} 
&\displaystyle{
(H_e)_{l,m}=\delta_{l,m}+\sum_{1\leq j\leq n} \,\Big(H^j_{l,m}\,z_j+\overline{H}^j_{m,l}\,\bar{z}_j\Big)}&
\\
\\
&\displaystyle{
+\sum_{1\leq j,k\leq n}\,\Big(H^{j,k}_{l,m}\,z_jz_k+\overline{H}^{j,k}_{m,l}\,\bar{z}_j\bar{z}_k+\hat{H}^{j,\bar{k}}_{l,m}\,z_j\bar{z}_k\Big)+O(|z|^3) }&  
\end{eqnarray*} 
les coefficients de la métrique hermitienne $h$ par rapport au repère $e$. La relation 
$$
A'_e=\overline{H}_e^{-1}(\partial_{_J }\overline{H}_e-\overline{A''_e}^{t}\,\overline{H}_e) $$
combinée avec les égalités $A'_e(x)=0, \,A''_e(x)=0$ implique alors $\bar{\partial}_{_{J}}H_e(x)=0$ et donc $H^j_{l,m}=0$ pour tout les indices $j,l,m$. Par rapport aux coordonnées choisies on a l'écriture
$$
(\partial_{_J } A''_e)_{m,l} =\sum_{1\leq j,k\leq n} \,(\partial_{_J } A''_e)^{j,\bar{k} }_{m,l}(0)\,dz_j\wedge d\bar{z}_k+O(|z|). 
$$
Considérons maintenant le changement de repère $\sigma =e \cdot g$ donné par la formule
$$
\sigma_l=e_l-\sum_{\substack{1\leq m\leq r_{\lambda}
\\
\\
1\leq j,k\leq n } }\,\Big(H^{j,k}_{l,m}\,z_jz_k+
(\partial_{_J }A''_e)^{j,\bar{k} }_{m,l}(0)\,z_j\bar{z}_k\Big)e _m\in{\cal E}(F^{\lambda }_{_{J}})(U_x).
$$
Un calcul élémentaire montre que les coefficients de la métrique hermitienne $h$ par rapport à ce repère s'écrivent sous la forme
$$
(H_{\sigma})_{l,m}=\delta _{l,m}+\sum_{j,k}\,H^{j,\bar{k}}_{l,m}\,z_j\bar{z}_k+O(|z|^3).
$$
Si $A''_{\sigma}$ désigne la forme de connexion relative au repère $\sigma$ on a la formule de changement de matrice de connexion $A''_{\sigma  }=g^{-1}(\bar{\partial}_{_J }g+A''_{e}\cdot g)$. Le fait que 
$A''_{e}(x)=0$ et $\bar{\partial}_{_J }g(x)=0$ implique alors l'égalité  $A''_{\sigma}(x)=0$.
De plus au point $x$ on a l'égalité 
$$
\partial_{_J } A''_{\sigma  }(x)=\partial_{_{J} } \bar{\partial}_{_J }g(x)+\partial_{_J }A''_e(x)=0.
$$
On déduit alors d'après la formule $\ref{ponctcurvat}$  que la courbure de Chern s'écrit au point $x$ sous la forme
$$
{\cal C}_h(F^{\lambda  }_{_{J}})_{|x}  =-\sum_{m,l}\,\partial_{_J }\bar{\partial}_{_J }h_{l,m}(x)\,\otimes\sigma^*_m\otimes_{_{J}}\sigma_l.      
$$
qui montre la validité de la formule $\eqref{curvat1}$. La formule $\eqref{curvat2}$ est une conséquence immédiate des identités  
$$
A'_{\sigma }(x) =\overline{H}_{\sigma} ^{-1}(\partial_{_J }\overline{H}_{\sigma } -\overline{A''_{\sigma} }^{t}\,\overline{H}_{\sigma} )(x)=0 
$$ et $\eqref{chernscalspec}$.\hfill$\Box$
\subsection{La courbure de Chern du fibré tangent d'une variété presque complexe}
Dans le cas du fibré tangent d'une variété presque complexe le tenseur de courbure de Chern $$
{\cal C}_{\omega }(T_{X,J}):=\Theta (D^{\omega }_{_{J} })^{1,1}\in {\cal E}(\Lambda ^{1,1}_{_{J}}T^*_X\otimes_{_{\C} } T^*_{X,J}\otimes_{_{\C} }T_{X,J})(X)  
$$
s'écrit sous la forme locale
\begin{eqnarray}\label{localcurvattg}
{\cal C}_{\omega }(T_{X,J})=\sum_{1\leq j,k,l,m \leq n}C^{j,k}_{l,m}\,\zeta ^*_j\wedge \bar{\zeta } ^*_k\otimes \zeta^*_m\otimes_{_{J}}\zeta _l.  
\end{eqnarray}
La notation $\alpha \otimes \zeta ^*_l\otimes_{_{J}}\zeta _m$ où $\alpha$ est une $(1,1)$-forme par exemple doit être interprétée sous la forme suivante. Si $\xi _1,\xi _2\in T_{X,x} \otimes_{_{\R}}\C$ et $\eta=\eta_l\zeta _l+\bar{\eta}_l\bar{\zeta }_l  \in T_{X,x} $  alors 
$$
\alpha \otimes \zeta ^*_l\otimes_{_{J}}\zeta _m(\xi_1,\xi _2,\eta)=\alpha(\xi _1,\xi _2)\,\eta_l\,\zeta _m
+\overline{\alpha(\xi _1,\xi _2)\,\eta_l\,\zeta _m}.
$$
En particulier la courbure de Chern du fibré tangent 
$$
{\cal C}^{\omega}_{_{{X,J} } }  \in {\cal E}(T_{X,J}^{*,\otimes 2}\otimes_{_{\C} }T_{X,-J}^{*,\otimes 2})(X)
$$ 
est définie par la formule 
$$
{\cal C}^{\omega}_{_{{X,J} } }(\xi _1\otimes\eta_1,\xi _2\otimes\eta_2)
:=h_{\omega}({\cal C}_{\omega }(T_{X,J})(\xi ^{1,0}_1,\xi ^{0,1}_2)\cdot \eta_1,\eta_2)
$$ 
pour tout champ de vecteurs réels $\xi _j,\eta_j\in{\cal E}(T_X)(U),\,j=1,2$, où $h_{\omega}$ est la forme hermitienne associée à $\omega$. On rappelle qu'elle est définie par la formule $h_{\omega}(\xi ,\eta):=\omega (\xi ,J\eta)-i\omega (\xi ,\eta)$. Le fait que ${\cal C}^{\omega}_{_{X,J} }$ soit une forme hermitienne sur le fibré $T^{\otimes 2}_{X,J}$ implique que la quantité ${\cal C}^{\omega}_{_{{X,J} } }(\xi \otimes\eta ,\xi \otimes\eta),\,\xi,\eta\in{\cal E}(T_X)(U)$ est réelle. On déduit alors les identités
$$
{\cal C}^{\omega}_{_{{X,J} } }(\xi \otimes\eta ,\xi \otimes\eta)=\omega ({\cal C}_{\omega }(T_{X,J})(\xi ,J\xi)\cdot \eta,\eta)
$$
et
$$
\omega ({\cal C}_{\omega }(T_{X,J})(\xi ,J\xi)\cdot \eta,J\eta)=0.
$$ 
La courbure de Chern du fibré tangent s'écrit en un point $x$ où le repère $(\zeta _k)_{k} \in{\cal E}(T^{1,0}_{X,J})^{\oplus n}(U)$ est choisie $\omega(x)$-orthonormé sous la forme
$$
{\cal C}^{\omega}_{_{{X,J} } }(x)=\sum_{1\leq j,k,l,m \leq n}C^{j,k}_{l,m}(x)\,\zeta ^*_j\otimes\zeta^* _m\otimes  \bar{\zeta } ^*_k\otimes \bar{\zeta }^*_l
$$
avec la relation de symétrie hermitienne $\overline{C^{j,k}_{l,m}(x)}=C^{k,j}_{m,l}(x)$. 
\\
\\
{\bf Remarque}. Le fait que la connexion de Chern soit hermitienne implique que en un point $x$ on a 
$\Theta (D^{\omega }_{_{J} })^{0,2}_{|x} =0$ si et seulement si $\Theta (D^{\omega }_{_{J} })^{2,0}_{|x} =0$. On peut montrer que $\Theta (D^{\omega }_{_{J} })^{0,2}_{|x} =0$ si le jet d'ordre un de la forme de torsion de la structure presque complexe est nul au point $x$.
\section{Coordonnées presque complexes d'ordre $N$ en un point} 
Soient $(z_1,...,z_n)$ des coordonnées locales $\ci$ centrées en $x\in X$ telles que le repère local $(\frac{\partial}{\partial z_1},...,\frac{\partial}{\partial z_n})$ soit une base complexe de $T^{1,0}_{X,J,x}$ au point $x$. On désigne par $M_J\in M_{2n,2n}({\cal E})$ la matrice de la structure presque complexe 
$J\in{\cal E}(End_{_{\C} } (T_X\otimes_{_{\R}}\C))(X)$ par rapport au repère complexe
$(\frac{\partial}{\partial z_1},...,\frac{\partial}{\partial z_n},\frac{\partial}{\partial \bar{z}_1},...,\frac{\partial}{\partial \bar{z}_n} )$. Le fait que $\overline{J}=J$ implique que la matrice $M_J$ s'écrit sous la forme: 
$$
M_J(z)=\left(
\begin{array}{cc} 
A(z)& \overline{B}(z)\\
B(z)&\overline{A}(z)  
\end{array}  
\right)
$$
On voit alors que la structure presque complexe s'exprime sous la forme:
$$
J(z)=\sum_{k,l}\Big(A_{k,l}(z)\,dz_l\otimes \frac{\partial}{\partial z_k} +B_{k,l}(z)\, dz_l\otimes\frac{\partial}{\partial \bar{z}_k}+\overline{B}_{k,l}(z)\, d\bar{z}_l\otimes \frac{\partial}{\partial z_k}
+\overline{A}_{k,l}(z)\, d\bar{z}_l\otimes\frac{\partial}{\partial \bar{z}_k}  \Big)
$$
avec $A(0)=iI_n,\;B(0)=0_n$. Si on suppose que la structure presque complexe est intégrable il existe d'après le théorème de Newlander-Nirenberg des coordonnées locales holomorphes $(z_1,...,z_n)$. La structure presque complexe s'écrit alors par rapport à ces coordonnées sous la forme
\begin{eqnarray}
J(z)=J_0=i\sum_k\Big(dz_k\otimes \frac{\partial}{\partial z_k}-d\bar{z}_k\otimes\frac{\partial}{\partial \bar{z}_k}\Big)\label{J0} 
\end{eqnarray}
autrement dit $A(z)\equiv i\,I_n,\;B(z)\equiv 0_n$.
Avec les notations introduites précédemment on a la proposition suivante.
\begin{prop}\label{JN} 
Pour tout point $x$ d'une variété presque complexe $(X,J)$ et pour tout entier $N\geq 2$ il existe des coordonnées $(z_1,...,z_n)$ de classe $\ci$ centrées en $x$ telles que  les matrices $A(z)$ et $B(z)$ de la structure presque complexe $J$ relatives à ces coordonnées admettent les développements asymptotiques
\begin{eqnarray}
&A (z)=i\,I_n+\displaystyle{\frac{i}{2}\sum_{\scriptstyle |\alpha +\beta |\leq N
\atop
\scriptstyle |\alpha |, |\beta |\geq 1} \,A^{\alpha ,\beta }  \,z^{\alpha} \bar{z}^{\beta } +O(|z|^{N+1} )}&\label{AN}
\\\nonumber
\\
&B(z)=\displaystyle{\sum_ {\scriptstyle |\alpha +\beta |\leq N
\atop
\scriptstyle |\alpha |\geq 1} }\,B^{\alpha ,\beta }  \,z^{\alpha} \bar{z}^{\beta } +O(|z|^{N+1} )&\label{BN}
\end{eqnarray}
où $A^{\alpha ,\beta } ,\,B^{\alpha ,\beta } \in M_{n,n}(\C)$ sont des matrices telles que les coefficients des matrices $B^{\alpha ,\beta }$ vérifient la propriété; $B^{\alpha ,\beta }_{k,l}  =0$ pour tout $l\geq \max\{k\in\{1,...,n\}\,|\,\alpha _k\not =0 \}$. Les matrices  $A^{\alpha ,\beta }$ sont obtenues à partir des matrices $B^{\alpha ,\beta }$, $($avec la convention $B^{0,\beta }:=0$$)$, grâce à la formule:
\begin{eqnarray}\label{formA}
\displaystyle{ A^{\alpha ,\beta }=\sum_{k=1}^{[|\alpha +\beta |/2]}\sum _{\scriptstyle \sum_{r=1}^k\,(\rho _r+\mu _r)=\alpha  
\atop
\scriptstyle  \sum_{r=1}^k\,(\lambda_r+\gamma_r )=\beta  } \,(-4)^{-(k-1)}\prod_{1\leq r\leq k }^{\longrightarrow}\overline{B}^{\lambda_r ,\mu _r }\cdot B^{\rho _r, \gamma_r}}
\end{eqnarray}
où le symbole $[c]$ désigne la partie entière de $c$ et le symbole de produit avec une flèche vers la droite désigne le produit non commutatif des termes qui sont écrits en ordre croissant de l'indice vers la droite. 
\end{prop}
(Remarquons que dans la formule $\eqref{formA}$ la convention $B^{0,\beta }=0$ implique que les sommes non nulles sont celles correspondantes aux multi-indices $|\lambda _r|,\,|\rho _r|\geq 1$).
\begin{defi}
Les coordonnées qui vérifient les propriétés de l'énonce de la proposition précédente seront appelées coordonnées presque complexes d'ordre $N$ en $x$ par rapport à la structure $J$.
\end{defi}  
Dans le cas particulier $N=3$ 
 la formule $\eqref{formA}$ s'écrit sous la forme;
$$
A^{\alpha ,\beta }= \sum_{\scriptstyle \mu +\rho =\alpha 
\atop
\scriptstyle \lambda +\gamma=\beta } \,\overline{B}^{\lambda }\cdot B^{\rho,\gamma}  
$$
On ré-énonce la proposition précédente dans le cas $N=3$ sous une forme plus explicite et pratique pour les calculs relatifs à la sous-section qui suivra.
\begin{coro}\label{J3}
Pour tout point $x$ d'une variété presque complexe $(X,J)$ il existe des coordonnées $(z_1,...,z_n)$ de classe $\ci$ centrées en $x$ telles que les matrices $A(z)$ et $B(z)$ de la structure presque complexe $J$ relatives à ces coordonnes admettent les développements asymptotiques
\begin{eqnarray}
&B(z)=\displaystyle{\sum_r \,B^rz_r+\sum_{r,s}}\, \Big(B^{r,s}\,z_rz_s +B^{r,\bar{s}}\,z_r\bar{z}_s \Big)&\nonumber
\\\nonumber
\\
&\displaystyle{+\sum_{r,s,t}\,\Big(B^{r,s,t}\,z_rz_sz_t+B^{r,s,\bar{t} }\,z_rz_s\bar{z}_t+B^{r,\bar{s},\bar{t} }\,z_r\bar{z}_s\bar{z}_t\Big) +O(|z|^4)}&\label{B3}
\\\nonumber
\\
&A (z)=i\,I_n+\displaystyle{\frac{i}{2}\sum_{r,s} \,\overline{B}^{r} \cdot B^{s} \,z_s\bar{z}_r+
\frac{i}{4 } \sum_{r,s,t}\,\Big(\overline{B}^{t,\bar{r} }\cdot B^{s}+\overline{B}^{t,\bar{s}}\cdot B^r+2\overline{B}^t\cdot B^{r,s}\Big)\, z_rz_s\bar{z}_t  }&\nonumber
\\\nonumber
\\
&\displaystyle{+\frac{i}{4 } \sum_{r,s,t}\,\Big(\overline{B}^t\cdot B^{r,\bar{s}}+\overline{B}^s\cdot B^{r,\bar{t} }+2 \overline{B}^{s,t}\cdot B^r\Big)\, z_r\bar{z}_s\bar{z}_t +O(|z|^4)  }&
 \label{A3}
\end{eqnarray}
où $B^r,\,B^{r,s},\,B^{r,\bar{s} },\,B^{r,s,t},\,B^{r,s,\bar{t}},\,B^{r,\bar{s},\bar{t} }\in M_{n,n}(\C)$ sont des matrices telles que $B^{r,s}$ soit symétrique par rapport aux indices $r,s$, $B^{r,s,t}$ par rapport à $r,s,t$, 
$B^{r,s,\bar{t} }$ par rapport à $r,s$, $B^{r,\bar{s},\bar{t}}$ par rapport à $s,t$ et $B^r_{k,l}=0$ pour $r\leq l$, $B^{r,s}_{k,l}=0 $ pour $r,s\leq l$, $B^{r,\bar{s}}_{k,l}=0$ pour $r\leq l$, $B^{r,s,t}_{k,l}=0$ pour 
$r,s,t\leq l$, $B^{r,s,\bar{t}}_{k,l}=0$ pour $r,s\leq l$, et $B^{r,\bar{s} ,\bar{t} }_{k,l}=0$ pour $r\leq l$. De plus si on considère l'expression locale de la forme de torsion de la structure presque complexe
\begin{eqnarray*}
&\displaystyle{\tau_{_{J}}=\sum_{1\leq k<l\leq n}[\zeta _k,\zeta _l]^{0,1}_{_{J}}\otimes \zeta ^*_k\wedge\zeta ^*_l 
=\sum_{\substack{1\leq k<l\leq n
\\
1\leq r\leq n }}\overline{N}^r_{k,l}\, \zeta ^*_k\wedge\zeta ^*_l \otimes \bar{\zeta}_r }&
\end{eqnarray*}
où $\zeta _l:=(\partial /\partial z_l)^{1,0}_{_{J} }\in {\cal E}(T^{1,0}_{X,J})(U_x), \,l=1,...,n$ est le repère locale du fibré des $(1,0)$-vecteurs $T^{1,0}_{X,J}$ issue des coordonnées $(z_1,...,z_n)$ on a l'expression
$$
\overline{N}^r_{k,l}(z)=\frac{i}{2}\,B^l_{r,k}+\frac{i}{2}\sum_s\Big[2(B^{l,s}_{r,k}-B^{k,s}_{r,l})\,z_s+ B^{l,\bar{s}}_{r,k} \,\bar{z}_s\Big]+O(|z|^2)
$$
pour tout $k<l$. Le jet d'ordre $k=0,1$ de la forme de torsion de la structure presque complexe au point $x$ est nul si et seulement si les coefficients $B_{*,*}(z)$ de la structure presque complexe relatifs aux coordonnées en question s'annulent à l'ordre $k+1$.
\end{coro}
Les coordonnées précédentes seront appelées coordonnées presque complexes d'ordre 3 au point $x$.
\\
\\
{\bf Preuve de la proposition \ref{JN}} 
\\
{\bf I) Les changements de coordonnées} 
\\ 
La condition $J^2=-\I$ est exprimée par les conditions locales $A^2=-I_n-\overline{B}\cdot B$ et $\overline{A}\cdot B=-B\cdot A$.
Le choix fait sur les coordonnées locales implique que relativement à celles-ci on a $J(0)=J_0,\;A(0)=i\,I_n,\;B(0)=0_n$. 
La relation $A^2=-I_n-\overline{B}\cdot B$ implique alors que la matrice $A$ admet un développement asymptotique du type $A(z)=i\,I_n+O(|z|^2)$.
Si $Z=\Phi (z)$ est un changement de coordonnées alors la matrice de la structure presque complexe
$$
{\cal M}_J(Z)=\left(
\begin{array}{cc} 
{\cal A} (Z)& \overline{\cal B }(Z)
\\
{\cal B} (Z)&\overline{\cal A}(Z)  
\end{array}  
\right)
$$
par rapport aux nouvelles coordonnées est donné par la formule ${\cal M}_J(Z)=d\Phi \cdot M_{J}(z)\cdot d\Phi^{-1}$. De manière explicite on a alors les formules
\begin{eqnarray}
&\displaystyle{
{\cal A}_{k,l}(Z):=\sum_{s,t}\,\Big(A_{s,t}(Z)\frac{\partial z_t}{\partial Z_l}\frac{\partial Z_k}{\partial z_s}
+B_{s,t}(Z)\frac{\partial z_t}{\partial Z_l}\frac{\partial Z_k}{\partial\bar{ z}_s}}&
\nonumber
\\\nonumber
\\
&\displaystyle{
+\overline{B}_{s,t}(Z)\frac{\partial \bar{ z}_t}{\partial Z_l}\frac{\partial Z_k}{\partial z_s}
+\overline{A}_{s,t}(Z)\frac{\partial \bar{ z}_t}{\partial Z_l}\frac{\partial Z_k}{\partial \bar{z}_s} \Big)}& \label{A} 
\end{eqnarray}
\begin{eqnarray}
&\displaystyle{
{\cal B}_{k,l}(z):=\sum_{s,t}\, \Big(A_{s,t}(Z)\frac{\partial z_t}{\partial Z_l}\frac{\partial \bar{Z} _k}{\partial z_s}
+B_{s,t}(Z)\frac{\partial z_t}{\partial z_l}\frac{\partial \bar{Z} _k}{\partial\bar{ z}_s}}&\nonumber
\\\nonumber
\\
&\displaystyle{
+\overline{B}_{s,t}(Z)\frac{\partial \bar{ z}_t}{\partial Z_l}\frac{\partial \bar{Z} _k}{\partial z_s}
+\overline{A}_{s,t}(Z)\frac{\partial \bar{ z}_t}{\partial Z_l}\frac{\partial \bar{Z} _k}{\partial \bar{z}_s} \Big).}&\label{B} 
\end{eqnarray}
Considérons maintenant  pour tout entier $N\geq 1$ les changements de coordonnées $Z=\Phi _N(z)$ 
$$
\displaystyle{ Z_k=z_k-\sum_ {\scriptstyle |\alpha +\beta | =N+1
\atop
\scriptstyle |\alpha |\geq 1} }\,\frac{i\, \overline{B}^{\alpha -\delta _{l(\alpha )},\beta  }_{k,l(\alpha )} }{2\alpha _{l(\alpha )} } \,z^{\beta} \bar{z}^{\alpha }   
$$
où  $l(\alpha ):=\max\{r\in\{1,...,n\}\,|\,\alpha _r\not =0 \}$ et les coefficients $B^{\alpha ,\beta },\,|\alpha +\beta |=N$ seront définis dans la suite. On considère aussi les changements inverses
$$
\displaystyle{ z_k=Z_k+\sum_ {\scriptstyle |\alpha +\beta | =N+1
\atop
\scriptstyle |\alpha |\geq 1} }\,\frac{ i\,\overline{B}^{\alpha -\delta _{l(\alpha )},\beta  }_{k,l(\alpha )} }{2\alpha _{l(\alpha )} } \,Z^{\beta} \bar{Z}^{\alpha } +O(|Z|^{2N+1} ) 
$$
On définit aussi
$$
{\cal B}^{\alpha-\delta _l ,\beta }_{k,l} :=B^{\alpha-\delta _l ,\beta }_{k,l}- \displaystyle{\alpha_l\cdot \frac{ B^{\alpha -\delta _{l(\alpha )},\beta  }_{k,l(\alpha )} }{\alpha _{l(\alpha )} }  }
$$
pour tout les multi-indices $\alpha $ tels que $\alpha _l\geq 1$. Avec la convention $0=\max \emptyset$, on a alors ${\cal B}^{\alpha,\beta }_{k,l}=0$ pour tout les multi-indices $|\alpha +\beta |= N$ tels que $l(\alpha )\leq l$. Avec la convention précédente on a en particulier ${\cal B}^{0,\beta }=0$ lorsque $|\beta |=N$.
On a les expressions suivantes pour les dérivées partielles:
\begin{eqnarray*}
&\displaystyle{\frac{\partial z_t}{\partial Z_l}=\delta _{t,l}+\sum_ {\scriptstyle |\alpha +\beta | =N+1
\atop
\scriptstyle |\alpha |,\,\beta _l\geq 1} \,\beta _l\cdot \frac{ i\,\overline{B}^{\alpha -\delta _{l(\alpha )},\beta  }_{t,l(\alpha )} }{2\alpha _{l(\alpha )} } \,Z^{\beta-\delta _l} \bar{Z}^{\alpha }   +O(|Z|^{2N} ) } &
\end{eqnarray*}
\begin{eqnarray*}
&\displaystyle{\frac{\partial z_t}{\partial\bar{Z}_l}=\sum_ {\scriptstyle |\alpha +\beta | =N+1
\atop
\scriptstyle \alpha _l\geq 1} }\,\alpha_l\cdot \frac{ i\,\overline{B}^{\alpha -\delta _{l(\alpha )},\beta  }_{t,l(\alpha )} }{2\alpha _{l(\alpha )} } \,Z^{\beta} \bar{Z}^{\alpha-\delta _l } +O(|Z|^{2N} )  &
\\
\\
&\displaystyle{\frac{\partial Z_k}{\partial z_s}=\delta _{s,k}-\sum_ {\scriptstyle |\alpha +\beta | =N+1
\atop
\scriptstyle |\alpha |,\,\beta _s\geq 1} \,\beta _s\cdot \frac{ i\,\overline{B}^{\alpha -\delta _{l(\alpha )},\beta  }_{k,l(\alpha )} }{2\alpha _{l(\alpha )} } \,Z^{\beta-\delta _s} \bar{Z}^{\alpha } +O(|Z|^{2N} ) }&
\\
\\
&\displaystyle{ \frac{\partial Z_k}{\partial \bar{z}_s}=-\sum_ {\scriptstyle |\alpha +\beta | =N+1
\atop
\scriptstyle \alpha _s\geq 1} }\,\alpha_s\cdot \frac{ i\,\overline{B}^{\alpha -\delta _{l(\alpha )},\beta  }_{k,l(\alpha )} }{2\alpha _{l(\alpha )} } \,Z^{\beta} \bar{Z}^{\alpha-\delta _s } +O(|Z|^{2N} )  &
\end{eqnarray*}
Nous allons montrer maintenant à l'aide d'une récurrence sur $N$, l'existence de coordonnées pour lesquelles les matrices $A(z)$ et $B(z)$ admettent les  développements asymptotiques $\eqref{AN} $ et $\eqref{BN} $ avec les conditions sur les coefficients $B^{\alpha ,\beta }_{k,l}$ expliquées dans l'énonce du lemme. On commence par effectuer le changement de coordonnées $Z=\Phi_1(z)$ où les matrices $B^{\alpha ,\beta },\,|\alpha +\beta |=1$ qui apparaissent dans la définition de tel changement sont celles du développement:
$$
B(z)=\sum_{|\alpha +\beta |=1} \,B^{\alpha ,\beta }\,z^{\alpha }\bar{z}^{\beta }+O(|z|^2)    
$$ 
(rappelons que $B(0)=0_n$). En substituant les expressions des dérivées partielles relatives au changement de coordonnées $Z=\Phi_1(z)$ et en tenant compte des développements asymptotiques des matrices $A(z)$ et $B(z)$ obtenues précédemment dans les expressions $\eqref{A}$, $\eqref{B}$ on aura, relativement aux nouvelles coordonnées, les développements asymptotiques suivants:
\begin{eqnarray*}
&\displaystyle{  {\cal A}_{k,l} (Z)=\sum_{s,t}\,A_{s,t}(Z)\frac{\partial z_t}{\partial Z_l}\frac{\partial Z_k}{\partial z_s}+O(|Z|^2)=i\,\delta _{k,l}+O(|Z|^2) }&
\\
\\
&\displaystyle{ {\cal B}_{k,l}(Z)=\sum_s\, i\,\Big(\frac{\partial z_s}{\partial Z_l}\frac{\partial \bar{Z}_k}{\partial z_s}-
\frac{\partial \bar{z} _s}{\partial Z_l}\frac{\partial \bar{Z}_k}{\partial \bar{z} _s}\Big)+B_{k,l}(Z)+O(|Z|^{2})= }&
\\
\\
&\displaystyle{=-\sum_ {\scriptstyle |\alpha +\beta | =2
\atop
\scriptstyle \alpha_l \geq 1} \,\alpha_l\cdot \frac{ B^{\alpha -\delta _{l(\alpha )},\beta  }_{k,l(\alpha )} }{\alpha _{l(\alpha )} } \,Z^{\alpha-\delta _l } \bar{Z}^{\beta}  +B_{k,l}(Z)+O(|Z|^{2})=}&
\\
\\
&\displaystyle{=\sum_{\scriptstyle |\alpha +\beta |= 2
\atop
\scriptstyle \alpha _l\geq 1} \,{\cal B}^{\alpha-\delta _l ,\beta }_{k,l}  \,Z^{\alpha-\delta _l} \bar{Z}^{\beta }  +O(|Z|^{2} ) }&
\end{eqnarray*}
Pour simplifier les notations dans les calculs qui suivront on va noter à partir de maintenant $A$ à la place de ${\cal A}$, $B$  à la place de ${\cal B}$ et $z$  à la place de $Z$. Avec ces notations on a alors que la matrice $B(z)$ peut être écrite sous la forme asymptotique $\eqref{BN}$ avec $N=1$ et les conditions correspondantes sur les coefficients $B^{\alpha ,\beta }_{k,l}$. La relation $A^2=-I_n-\overline{B}\cdot B$ entraîne alors que la matrice $A(z)$ admet le développement asymptotique $\eqref{AN}$ avec $N=2$. Supposons maintenant qu'il existe des coordonnées telles que la matrice $B(z)$ admette le développement $\eqref{BN}$ relativement à l'entier $N-1,\,N\geq 2$. On peut alors écrire le développement asymptotique suivant:
$$
B(z)=\displaystyle{\sum_ {\scriptstyle |\alpha +\beta |\leq N-1
\atop
\scriptstyle |\alpha |\geq 1} }\,B^{\alpha ,\beta }  \,z^{\alpha} \bar{z}^{\beta }+\sum_{ |\alpha +\beta |=N}\,B^{\alpha ,\beta }  \,z^{\alpha} \bar{z}^{\beta }  +O(|z|^{N+1} )
$$
relativement aux coordonnées en question, où $B^{\alpha ,\beta }_{k,l}=0$ pour $l\geq l(\alpha ),\,|\alpha +\beta |\leq N-1,\,|\alpha |\geq 1$. L'expression précédente de $B(z)$ combinée avec la relation 
$A^2=-I_n-\overline{B}\cdot B$ implique que la matrice $A(z)$ s'écrit sous la forme $\eqref{AN}$. On considère maintenant le changement de coordonnées $Z=\Phi_N(z)$ où les matrices 
$B^{\alpha ,\beta } ,\,|\alpha +\beta |=N$ qui apparaissent dans la définition de tel changement sont celles qui apparaissent dans l'expression asymptotique précédente de $B(z)$. 
Par rapport aux nouvelles coordonnées les matrices $A(z)$ et $B(z)$ admettent les développements asymptotiques suivants:
\begin{eqnarray*}
&\displaystyle{ {\cal A}_{k,l} (Z)=\sum_{s,t}\,A_{s,t}(Z)\frac{\partial z_t}{\partial Z_l}\frac{\partial Z_k}{\partial z_s}+O(|Z|^{N+1})=}&
\\
\\
&\displaystyle{=i\delta _{k,l}+\frac{i}{2}\sum_{\scriptstyle |\alpha +\beta |\leq N
\atop
\scriptstyle |\alpha |, |\beta |\geq 1} \,A^{\alpha ,\beta }_{k,l}   \,Z^{\alpha} \bar{Z}^{\beta } +O(|Z|^{N+1} )  } ,&
\\
\\
&\displaystyle{ {\cal B}_{k,l}(\zeta)=\sum_s\, i\,\Big(\frac{\partial z_s}{\partial Z_l}\frac{\partial \bar{Z}_k}{\partial z_s}-
\frac{\partial \bar{z} _s}{\partial Z_l}\frac{\partial \bar{Z}_k}{\partial \bar{z} _s}\Big)+B_{k,l}(Z)+O(|Z|^{N+1})= }&
\\
\\
&\displaystyle{=-\sum_ {\scriptstyle |\alpha +\beta | =N+1
\atop
\scriptstyle \alpha_l \geq 1} \,\alpha_l\cdot \frac{ B^{\alpha -\delta _{l(\alpha )},\beta  }_{k,l(\alpha )} }{\alpha _{l(\alpha )} } \,Z^{\alpha-\delta _l } \bar{Z}^{\beta}  +B_{k,l}(Z)+O(|Z|^{N+1})=}&
\end{eqnarray*}
\begin{eqnarray*}
&\displaystyle{=\sum_ {\scriptstyle |\alpha +\beta |\leq N-1
\atop
\scriptstyle |\alpha |\geq 1} \,B^{\alpha ,\beta }_{k,l}   \,Z^{\alpha} \bar{Z}^{\beta }+\sum_{\scriptstyle |\alpha +\beta |= N+1
\atop
\scriptstyle \alpha _l\geq 1} \,{\cal B}^{\alpha-\delta _l ,\beta }_{k,l}  \,Z^{\alpha-\delta _l} \bar{Z}^{\beta }  +O(|Z|^{N+1} ) }.&
\end{eqnarray*}
De la même façon que précédemment, on va noter à partir de maintenant $A$ à la place de ${\cal A}$, $B$  à la place de ${\cal B}$ et $z$  à la place de $Z$. Avec ces notations on obtient en conclusion que les matrices $A(z)$ et $B(z)$ peuvent être écrites sous les formes asymptotiques $\eqref{AN}$ et $\eqref{BN}$, avec les conditions correspondantes sur les coefficients $B^{\alpha ,\beta }_{k,l}$. 
\\
\\
{\bf II) Preuve de la formule $\eqref{formA}$} 
\\
On montre maintenant la formule  $\eqref{formA}$ à l'aide d'une récurrence sur $N\geq 2$. Pour simplifier les notations dans les calculs qui suivront on utilisera les conventions $A^{\alpha ,0}=A^{0,\beta }=0$.  En tenant compte des expressions $\eqref{AN}$ et $\eqref{BN}$  pour $2\leq N\leq 3$ on peut écrire la relation $A^2=-I_n-\overline{B}\cdot B$ sous la forme
\begin{eqnarray*}
&\displaystyle{-I_n-\sum_{|\alpha +\beta |\leq N}\,A^{\alpha ,\beta }\,z^{\alpha }\bar{z}^{\beta}+O(|z|^{N+1})=}&
\\
\\
&\displaystyle{=-I_n-\sum_{|\alpha +\beta |\leq N}\,\Big(\sum_{\scriptstyle \mu +\rho =\alpha 
\atop
\scriptstyle \lambda +\gamma=\beta } \,\overline{B}^{\lambda }\cdot B^{\rho ,\gamma}    \Big)\,z^{\alpha }\bar{z}^{\beta}+O(|z|^{N+1})  } ,&    
\end{eqnarray*}
(rappelons qu'on utilise la convention $B^{0,\beta }=0$). On a alors 
$$
A^{\alpha ,\beta }= \sum_{\scriptstyle \mu +\rho =\alpha 
\atop
\scriptstyle \lambda +\gamma=\beta } \,\overline{B}^{\lambda }\cdot B^{\rho ,\gamma}  
$$
pour $2\leq |\alpha +\beta |\leq 3$, qui n'est rien d'autre que la formule $\eqref{formA}$ dans les cas particuliers en considération. Nous supposons maintenant avoir montré la formule $\eqref{formA}$ pour $2\leq |\alpha +\beta |\leq N$. Comme précédemment la relation $A^2=-I_n-\overline{B}\cdot B$ s'écrit, à l'aide des expressions $\eqref{AN}$ et $\eqref{BN}$ pour $N+1$, sous la forme:
\begin{eqnarray*}
&\displaystyle{-I_n-\sum_{|\alpha +\beta |\leq N+1}\,A^{\alpha ,\beta }\,z^{\alpha }\bar{z}^{\beta}-\frac{1}{4}\sum_{|\alpha +\beta |\leq N+1}\,\Big(\sum_{\scriptstyle \lambda  +\rho =\alpha 
\atop
\scriptstyle \mu  +\gamma=\beta }\,A^{\lambda }\cdot A^{\rho ,\gamma}\Big)\,z^{\alpha }\bar{z}^{\beta}+  O(|z|^{N+2})=}&
\end{eqnarray*}
\begin{eqnarray*}
&\displaystyle{=-I_n-\sum_{|\alpha +\beta |\leq N+1}\,\Big(\sum_{\scriptstyle \mu +\rho =\alpha 
\atop
\scriptstyle \lambda +\gamma=\beta } \,\overline{B}^{\lambda }\cdot B^{\rho ,\gamma}    \Big)\,z^{\alpha }\bar{z}^{\beta}+O(|z|^{N+2})  }&    
\end{eqnarray*}
Cette identité implique que pour tout $\alpha ,\,\beta ,\,|\alpha +\beta |= N+1$ on a:
$$
A^{\alpha ,\beta }=\sum_{\scriptstyle \mu +\rho =\alpha 
\atop
\scriptstyle \lambda +\gamma=\beta } \,\overline{B}^{\lambda }\cdot B^{\rho ,\gamma}-\frac{1}{4}  \sum_{\scriptstyle \lambda  +\rho =\alpha 
\atop
\scriptstyle \mu  +\gamma=\beta }\,A^{\lambda }\cdot A^{\rho ,\gamma}  .
$$
En rappelant la convention $A^{0,\beta }=A^{\alpha ,0}=0$ on a que les termes non nuls de la dernière somme sont les termes relatifs aux multi-indices $|\lambda +\mu |,|\rho +\gamma|\leq N$. En utilisant l'hypothèse récursive relativement à l'expression  $\eqref{formA}$ on peut écrire l'expression précédente de la matrice $A^{\alpha ,\beta }$ sous la forme;
\begin{eqnarray*}
A^{\alpha ,\beta }=\sum_{\scriptstyle \mu +\rho =\alpha 
\atop
\scriptstyle \lambda +\gamma=\beta } \,\overline{B}^{\lambda} B^{\rho ,\gamma}-\qquad\qquad\qquad\qquad\qquad\qquad\qquad\qquad
\\
\\
-\frac{1}{4}\sum_{\substack{
\lambda  +\rho =\alpha 
\\
\mu  +\gamma=\beta
\\ 
\\
1\leq k_1\leq [|\lambda +\mu |/2]
\\ 
\\
 1\leq k_2\leq [|\rho +\gamma|/2]} }\; \sum _{\substack{ \sum_{r_1=1}^{k_1} \,(\rho _{r_1} +\mu _{r_1})=\lambda   
\\
\\
\sum_{r_1=1}^{k_1} \,(\lambda_{r_1}+\gamma_{r_1} )=\mu
\\
\\
\sum_{r_2=1}^{k_2} \,(\rho _{r_2}+\mu _{r_2})=\rho 
\\
\\
\sum_{r_2=1}^{k_2} \,(\lambda_{r_2}+\gamma_{r_2} )=\gamma }   } \,(-4)^{-(k_1+k_2-2)}
\prod_{1\leq r_1\leq k_1 }^{\longrightarrow}\overline{B}^{\lambda_{r_1} ,\mu _{r_1} }
 B^{\rho _{r_1}, \gamma_{r_1}} 
\prod_{1\leq r_2\leq k_2 }^{\longrightarrow}\overline{B}^{\lambda_{r_2} ,\mu _{r_2} }
 B^{\rho _{r_2}, \gamma_{r_2}}.
\end{eqnarray*}
En analysant l'ensemble des indices qui apparaissent sous les sommes précédentes on s'aperçoit de la validité de l'expression $\eqref{formA}$ relativement aux multi-indices $\alpha ,\,\beta $ en considération. \hfill $\Box$
\\
\\
{\bf Preuve du corollaire \ref{J3} }
\\
Le repère local $\zeta _k=(\partial /\partial z_k)^{1,0}_{_{J} }, \,k=1,...,n$ s'écrit sous la forme
\begin{eqnarray*}
&\displaystyle{
\zeta _k=\frac{1}{2}\frac{\partial}{\partial z_k}   -\frac{i}{2}\sum_r\Big(A_{r,k}\frac{\partial}{\partial z_r}  +
B_{r,k}\frac{\partial}{\partial \bar{z}_r}\Big)=}&
\\
\\
&\displaystyle{
=\frac{\partial}{\partial z_k} +\frac{1}{4}\sum_{p,h,t,j}\,\overline{B}^h_{t,j} B^p_{j,k}\,z_p\bar{z}_h  \frac{\partial}{\partial z_t} -\frac{i}{2}\sum_t\, {\bf jet}_2B_{t,k}(z)\,\frac{\partial}{\partial \bar{z}_t}  +O(|z|^3)   }&
\end{eqnarray*}
où ${\bf jet}_2\, B_{t,l}(z)$ désigne le jet d'ordre $2$ du coefficient $B_{t,l}$ de la structure presque complexe $J$ par rapport aux coordonnées en question.
On déduit alors facilement l'expression suivante pour le crochet
$$
[\zeta _k,\zeta _l]=\frac{i}{2}\sum_r\,\Big\{B^l_{r,k}+\sum_p\Big[2(B^{l,p}_{r,k}-B^{k,p}_{r,k})\,z_p+B^{l,\bar{p}}_{r,k}\,\bar{z}_p\Big] \Big\} \frac{\partial}{\partial \bar{z}_r} -\frac{1}{4}\sum_{pr,j}\,\overline{B}^p_{r,j}B^l_{j,k}\,\bar{z}_p\frac{\partial}{\partial z_r}+O(|z|^2).        
$$
En tenant compte de l'expression de la structure presque complexe à l'ordre un
$$
J(z)=i\sum_k\,\Big(dz_k\otimes\frac{\partial}{\partial z_k}-d\bar{z}_k\otimes\frac{\partial}{\partial \bar{z}_k}\Big)+
\sum_{k,l,p}\Big(B^p_{k,l}\,z_p\,dz_l\otimes \frac{\partial}{\partial \bar{z}_k}+
\overline{B}^p_{k,l}\,\bar{z}_p\,d\bar{z}_l\otimes \frac{\partial}{\partial z_k} \Big)+O(|z|^2) 
$$
on obtient l'expression
$$
J[\zeta _k,\zeta _l]=\frac{1}{2}\sum_r\,\Big\{B^l_{r,k}+\sum_p\Big[2(B^{l,p}_{r,k}-B^{k,p}_{r,k})\,z_p+B^{l,\bar{p}}_{r,k}\,\bar{z}_p\Big] \Big\} \frac{\partial}{\partial \bar{z}_r} +\frac{i}{4}\sum_{pr,j}\,\overline{B}^p_{r,j}B^l_{j,k}\,\bar{z}_p\frac{\partial}{\partial z_r}+O(|z|^2).
$$
On déduit alors l'expression
$$
[\zeta _k,\zeta _l]^{0,1}_{_{J} }  =\frac{i}{2}\sum_r\,\Big\{B^l_{r,k}+\sum_p\Big[2(B^{l,p}_{r,k}-B^{k,p}_{r,k})\,z_p+B^{l,\bar{p}}_{r,k}\,\bar{z}_p\Big] \Big\} \frac{\partial}{\partial \bar{z}_r} -\frac{1}{4}\sum_{pr,j}\,\overline{B}^p_{r,j}B^l_{j,k}\,\bar{z}_p\frac{\partial}{\partial z_r}+O(|z|^2).        
$$
En tenant compte de l'expression
$$
\bar{\zeta}_r= \frac{\partial}{\partial \bar{z}_r} +\frac{i}{2}\sum_{s,p}\,\overline{B}^p_{s,r}\,\bar{z}_p\frac{\partial}{\partial z_s}+   O(|z|^2)
$$
on déduit l'expression voulue pour les coefficients $\overline{N}^*_{*,*}$ de la forme de torsion de la structure presque complexe. Ces coefficients s'annulent à l'ordre $k=0,1$ si et seulement si les coefficients $B_{*,*}(z)$ de la structure presque complexe s'annulent à l'ordre $k+1$. En effet supposons que $T^{k,l,s}_r:= B^{l,s}_{r,k}-B^{k,s}_{r,l}$ soit nul pour tout les indices $k,l,s,r$. Si $k$ ou $l$ est le maximum de l'ensemble $\{k,l,s\}$ alors on a immédiatement 
$B^{l,s}_{r,k}=B^{k,s}_{r,l}=0$. Sinon, $s=\max\{k,l,s\}$ et donc $T^{k,s,l}_r=B^{s,l}_{r,k}=B^{l,s}_{r,k}=0$. \hfill $\Box$
\\
\\
Le calcul fait dans la preuve du corollaire $\ref{J3}$ montre que $M^*=O(|z|^2)$. Une conséquence immédiate des formules $\eqref{A} $ et $\eqref{B}$  est le corollaire suivant.
\begin{coro}\label{corholo}
Soient $(z_1,...,z_n)$ des coordonnées presque complexes à l'ordre $N\geq 1$ en un point $x$ et soit 
$Z_k=z_k+\sum_{|\alpha |=N+1}\,C^k_{\alpha }\,z^{\alpha}$ un changement de coordonnées holomorphe. Alors les coordonnées $(Z_1,...,Z_n)$ sont presque complexes à l'ordre $N$ en $x$ et les coefficients $B^{*,*}_{*,*}$ du jet d'ordre $N$ de la structure presque complexe par rapport aux nouvelles coordonnées sont les mêmes que les coefficients relatifs aux coordonnées $(z_1,...,z_n)$. 
\end{coro}   
\section{Expression asymptotique normale à l'ordre un d'une connexion de Chern sur le fibré tangent}\label{Exp-asymp-Chern} 
Le lemme suivant est nécessaire pour le calcul asymptotique du flot géodésique induit par une connexion de Chern sur le fibré tangent. L'expression asymptotique du flot de Chern est utile pour une technique de régularisation globale des $(1,1)$-courants positifs du type $i\partial_{_{J}}\bar{\partial}_{_{J}}u$ sur les variétés presque complexes (voir le chapitre trois pour plus de détails). Ce lemme et celui qui suivra montrent de façon optimale combien on est loin du cas Kählerien, où on dispose de coordonnées géodésiques complexes centrées en un point.
\begin{lem}\label{Exp-asym-Conn-Chern} 
Soit $(X,J)$ une variété presque complexe, $\omega\in {\cal E}(\Lambda ^{1,1}_{_J}T_X^*)(X)$ une métrique hermitienne et soient $(z_1,...,z_n)$  des coordonnées presque complexes d'ordre $N\geq 2$ en un point $x$ telles que le repère normal $(\zeta _k+\bar{\zeta}_k)_k\in{\cal E}(T_X)(U_x),\,\zeta _k=(\partial/\partial z_k)^{1,0}_{_{J} }$ soit $\omega (x)$-orthonormé. La métrique $\omega $ s'écrit alors sous la forme 
\begin{eqnarray}\label{jetmetri}  
&\displaystyle{
\omega =\frac{i}{2}\sum_{l,m}\,\Big[h_{l,m}+
\frac{i}{4}\sum_{j,k,r} \,B^j_{r,l}\overline{B}^k_{r,m}\,z_j\bar{z}_k  \Big]\,dz_l\wedge d\bar{z}_m }&\nonumber
\\\nonumber
\\  
&\displaystyle{-\frac{1}{4}\sum_{l,m}\,{\bf jet}_2 B_{l,m}(z)\,dz_l\wedge dz_m
-\frac{1}{4}\sum_{l,m} \,\overline{{\bf jet}_2 B_{l,m}(z)}\,d\bar{z}_l\wedge d\bar{z}_m+O(|z|^3)},&
\end{eqnarray}
où
$$
h_{l,m}=\delta _{l,m}+\sum_p\,\Big(H^p_{l,m}\,z_p+\overline{H}^p_{m,l}\,\bar{z}_p\Big)
+\sum_{p,h}\,\Big(H^{p,h}_{l,m}\,z_pz_h+\overline{H}^{p,h}_{m,l}\,\bar{z}_p\bar{z}_h+H^{p,\bar{h}}_{l,m}\,z_p\bar{z}_h\Big)+O(|z|^3) 
$$
et ${\bf jet}_2\, B_{l,m}$ désigne le jet d'ordre $2$ du coefficient $B_{l,m}$ de la structure presque complexe $J$ par rapport aux coordonnées en question. Pour tous champs de vecteurs réels
$
\eta= \sum_k\,(\eta _k \frac{\partial}{\partial z_k}+\bar{\eta}_k\frac{\partial}{\partial \bar{z}_k})
$
$
\in{\cal E}(T_X)(U_x)
$
on a l'expression asymptotique de la dérivée de Chern
$$
D^{\omega}_{_{J} }\,\eta=\sum_k\,\Big[d\eta_k+\sum_l\,\Big(E_{k,l}\,\eta_l-
\frac{i}{2}(d\;\overline{ {\bf jet}_2B_{k,l} }\,) \,\bar{\eta}_l\Big)\Big]\otimes_{_{J_0} }\frac{\partial}{\partial z_k}+O(|z|^2),
$$
où 
$$
E_{k,l}:=\sum_p\,\Big[H^p_{l,k}+\sum_h\,(S^{p,h}_{k,l}\,z_h+S^{p,\bar{h}}_{k,l}\,\bar{z}_h)\Big]\,dz_p+
\sum_{p,h}\,(S^{\bar{p},h}_{k,l}\,z_h+S^{\bar{p},\bar{h}}_{k,l}\,\bar{z}_h)\,d\bar{z}_p.
$$
Les coefficients $S^{*,*}_{k,l}$ sont données par les formules
\begin{align*}
&S^{\bar{p},h}_{k,l}=\frac{1}{4}\sum_j\,(\overline{B}^j_{k,p}-\overline{B}^p_{k,j})B^h_{j,l}  \,, 
&\quad 
S^{\bar{p},\bar{h}}_{k,l}=\frac{i}{2}\overline{B}^{h,\bar{l}}_{k,p}
-\frac{i}{2}\sum_j\,H^j_{l,k}\overline{B}^h_{j,p}\,,  
\\
&S^{p,h}_{k,l}=2H^{p,h}_{l,k}-\frac{i}{2}B^{h,\bar{k}}_{l,p}
-\sum_j\,\Big(H^p_{l,j}H^h_{j,k}+\frac{i}{2}\overline{H}^j_{k,l}B^h_{j,p}\Big)\,,
&\quad 
S^{p,\bar{h}}_{k,l}=-C^{p,h}_{k,l}(0)-\frac{1}{4}\sum_j\,\overline{B}^j_{k,h}B^p_{j,l} ,
\end{align*}
où $C^{p,h}_{k,l}(0)$ sont les coefficients de la courbure de Chern 
$$
{\cal C}_{\omega }(T_{X,J})=\sum_{j,k,m,l}\,C^{j,k}_{m,l}(0)\,dz_j\wedge d\bar{z}_k\otimes dz_l
\otimes_{_{J_0} }\frac{\partial}{\partial z_m}    +O(|z|).    
$$
au point $x$. Ils sont donnés par la formule
\begin{eqnarray}\label{valcurvat} 
C^{j,k}_{m,l}(0)=-H^{j,\bar{k}}_{l,m}+ \frac{1}{4}\sum_r\,\Big[4H^j_{l,r}\overline{H}^k_{m,r}+(\overline{B}^k_{m,r} -\overline{B}^r_{m,k})B^j_{r,l} 
+(B^j_{l,r}-B^r_{l,j})\overline{B}^k_{r,m}\Big].
\end{eqnarray}
\\
\end{lem}
$Preuve$
\\
On déduit facilement d'après la preuve du corollaire $\ref{J3}$ l'expression asymptotique à l'ordre deux du repère $(\zeta _l)_l$ et du repère dual $(\zeta ^*_l)_{l}$. On a les expressions asymptotiques suivantes.
\begin{eqnarray}
&\displaystyle{\zeta _l=\frac{\partial}{\partial z_l} +\frac{1}{4}\sum_{p,h,t,j}\,\overline{B}^h_{t,j} B^p_{j,l}\,z_p\bar{z}_h  \frac{\partial}{\partial z_t} -\frac{i}{2}\sum_t\, {\bf jet}_2B_{t,l}(z)\,\frac{\partial}{\partial \bar{z}_t}  +O(|z|^3)   }&\label{reptan} 
\\\nonumber
\\
&\displaystyle{\zeta^* _l=dz_l-\frac{i}{2}\sum_t\, \overline{ {\bf jet}_2B_{l,t}(z) } \,d\bar{z}_t+O(|z|^3) }.&\label{repdual} 
\end{eqnarray}
En tenant compte de cette dernière expression on déduit que la métrique
$
\omega =\frac{i}{2}\sum_{l,m}h_{l,m}\,\zeta _l^*\wedge\bar{\zeta}^*_m  
$
s'écrit sous la forme $\eqref{jetmetri}$.
On calcule maintenant les expressions asymptotiques des coefficients $U^*$, définis dans la section $1$, relativement au repère $\zeta _l=(\partial /\partial z_l)^{1,0}_{_{J} }, \,l=1,...,n$. Pour tout indice $k,h$ on a
\begin{eqnarray*}
&\displaystyle{[\zeta _k,\bar{\zeta} _h]=\sum_{r,l} \,\Big[\frac{i}{2}\overline{B}^{l,\bar{k}}_{r,h}\,\bar{z}_l 
+\frac{1}{4}\sum_j\,(\overline{B}^j_{r,h}-\overline{B}^h_{r,j})B^l_{j,k} \,z_l \Big]\frac{\partial}{\partial z_r }     }&
\\
\\
&\displaystyle{+\sum_{r,l} \,\Big[\frac{i}{2}B^{l,\bar{h}}_{r,k}\,z_l+\frac{1}{4}\sum_j\,(B^k_{r,j}-B^j_{r,k})\overline{B}^l_{j,h}  \,\bar{z}_l\Big]\frac{\partial}{\partial \bar{z}_r }+O(|z|^2)}& 
\end{eqnarray*}
En tenant compte de l'expression de la structure presque complexe à l'ordre un
$$
J(z)=i\sum_k\,\Big(dz_k\otimes\frac{\partial}{\partial z_k}-d\bar{z}_k\otimes\frac{\partial}{\partial \bar{z}_k}\Big)+
\sum_{k,l,p}\Big(B^p_{k,l}\,z_p\,dz_l\otimes \frac{\partial}{\partial \bar{z}_k}+
\overline{B}^p_{k,l}\,\bar{z}_p\,d\bar{z}_l\otimes \frac{\partial}{\partial z_k} \Big)+O(|z|^2) 
$$
on obtient l'expression
\begin{eqnarray*}
&\displaystyle{J[\zeta _k,\bar{\zeta} _h]=\sum_{r,l} \,\Big[-\frac{1}{2}\overline{B}^{l,\bar{k}}_{r,h}\,\bar{z}_l 
+\frac{i}{4}\sum_j\,(\overline{B}^j_{r,h}-\overline{B}^h_{r,j})B^l_{j,k} \,z_l \Big]\frac{\partial}{\partial z_r }     }&
\\
\\
&\displaystyle{+\sum_{r,l} \,\Big[\frac{1}{2}B^{l,\bar{h}}_{r,k}\,z_l-\frac{i}{4}\sum_j\,(B^k_{r,j}-B^j_{r,k})\overline{B}^l_{j,h}  \,\bar{z}_l\Big]\frac{\partial}{\partial \bar{z}_r }+O(|z|^2)}& 
\end{eqnarray*}
On a alors
$$
[\zeta _k,\bar{\zeta} _h]^{1,0}_{_{J} } =\sum_{r,l} \,\Big[\frac{1}{4}\sum_j\,(\overline{B}^j_{r,h}-\overline{B}^h_{r,j})B^l_{j,k} \,z_l+\frac{i}{2}\overline{B}^{l,\bar{k}}_{r,h}\,\bar{z}_l  \Big]\frac{\partial}{\partial z_r } +O(|z|^2).
$$
En tenant compte de l'expression asymptotique à l'ordre un du repère $(\zeta_k )_k$ on déduit l'expression
\begin{eqnarray}\label{conexionA} 
U^r_{k,h}(z)=\sum_l \,\Big[\frac{1}{4}\sum_j\,(\overline{B}^j_{r,h}-\overline{B}^h_{r,j})B^l_{j,k} \,z_l+\frac{i}{2}\overline{B}^{l,\bar{k}}_{r,h}\,\bar{z}_l  \Big] +O(|z|^2),
\end{eqnarray}
qui nous donne l'expression normale asymptotique à l'ordre un de la forme de connexion $A''_{\zeta }$ relative au repère normal $\zeta _k=(\partial/\partial z_k)^{1,0}_{_{J} }$. Nous calculons maintenant l'expression asymptotique à l'ordre un de la forme de connexion $A'_{\zeta }$ à l'aide de l'expression précédente de la forme $A''_{\zeta }$. La matrice inverse $H^{-1}=(h^{r,k})$ admet le developpement asymptotique suivant.$$
h^{r,k}=\delta_{r,k}-\sum_j\,(H^j_{r,k}\,z_j+\overline{H}^j_{k,r}\,\bar{z}_j)+O(|z|^2).
$$
En utilisant l'expression de la forme $A'_{\zeta }$ obtenue dans la preuve du théorème $\ref{teorconchern}$ on déduit l'expression
$$
(A'_{\zeta })_{k,l}=\sum_r\,h^{r,k}\partial_{_{J}}h_{l,r}+\sum_p\overline{U}^l_{k,p}\,dz_p+O(|z|^2),
$$
avec $\partial_{_{J}}h_{l,r}=\sum_p\,(\zeta_p\,.h_{l,r})\zeta^*_p$, où
$$
\zeta_p\,.h_{l,r}=H^p_{l,r}+\sum_h\,\Big[\Big(2H^{p,h}_{l,r}-
\frac{i}{2}\sum_t\,\overline{H}^t_{r,l}B^h_{t,p}\Big)\,z_h+H^{p,\bar{h}}_{l,r}\,\bar{z}_h\Big]+O(|z|^2).
$$
En utilisant l'expression du jet d'ordre un du repère $(\zeta^*_p)_p$ on obtient l'expression
\begin{eqnarray*}
&\displaystyle{
\partial_{_{J}}h_{l,r}=\sum_p\,\Big\{H^p_{l,r}+\sum_h\,\Big[\Big(2H^{p,h}_{l,r}-
\frac{i}{2}\sum_t\,\overline{H}^t_{r,l}B^h_{t,p}\Big)\,z_h+H^{p,\bar{h}}_{l,r}\,\bar{z}_h\Big]\Big\}\,dz_p}&\\
\\
&\displaystyle{
-\frac{i}{2}\sum_{p,t,h}\,H^t_{l,r}\overline{B}^h_{t,p}\,\bar{z}_h\,d\bar{z}_p+O(|z|^2).}&
\end{eqnarray*}
On déduit alors l'expression asymptotique à l'ordre un de la forme de connexion de Chern $A_{\zeta }=A'_{\zeta }+A''_{\zeta }$.
$$
A_{\zeta }=\partial_{_{J}}h_{l,k}-\sum_{p,r,j}\,H^p_{l,r}(H^j_{r,k}\,z_j+\overline{H}^j_{k,r}\,\bar{z}_j)\,dz_p+\sum_p\,(\overline{U}^l_{k,p}\,dz_p-U^k_{l,p}\,d\bar{z}_p).
$$
La matrice de la forme de connexion de l'extension 
$
D^{\omega}_{_{J} }:{\cal E}(T_X\otimes_{_{\R}}\C) \longrightarrow {\cal E}(T^*_X\otimes_{_{\R}}(T_X\otimes_{_{\R}}\C))
$
de la connexion de Chern au complexifié du fibré tangent $T_X\otimes_{_{\R}}\C$ par rapport au repère 
$(\zeta_k,\bar{\zeta}_k)_k$ est
$$
A_{\zeta,\bar{\zeta}}=\left(
\begin{array}{cc} 
A_{\zeta }& 0_n
\\
0_n& \overline{A}_{\zeta } 
\end{array}  
\right).
$$
On doit maintenant calculer la matrice $A_z$ de la forme de connexion de l'extension de la connexion de Chern par rapport au repère $(\frac{\partial}{\partial z_k},\frac{\partial}{\partial \bar{z}_k})_k$ du complexifié du fibré tangent $T_X\otimes_{_{\R}}\C$. La formule $\eqref{reptan}$ nous donne l'expression asymptotique de la matrice $g^{-1}$ du changement de repère $(\zeta_k,\bar{\zeta}_k)_k=(\frac{\partial}{\partial z_k},\frac{\partial}{\partial \bar{z}_k})_k\cdot g^{-1}$. Les expressions asymptotiques à l'ordre deux des matrices $g$ et $g^{-1}$ sont les suivantes 
$$
g=\left(
\begin{array}{cc} 
I_n                          & -\frac{i}{2}\overline{{\bf jet}_2 B}
\\
\\
\frac{i}{2}{\bf jet}_2 B & I_n
\end{array}  
\right)+O(|z|^3),
\qquad
g^{-1}=\left(
\begin{array}{cc} 
T                         & \frac{i}{2}\overline{{\bf jet}_2 B}
\\
\\
-\frac{i}{2}{\bf jet}_2 B & \overline{T}
\end{array}  
\right)+O(|z|^3),
$$
où $T_{k,l}=\delta_{k,l}+\frac{1}{4}\sum_{p,h,j}\,\overline{B}^h_{k,j} B^p_{j,l}\,z_p\bar{z}_h$. La matrice de la forme de connexion qu'on cherche est donnée par la formule $A_z=g^{-1}(dg+A_{\zeta,\bar{\zeta}}\,g)$. On a alors les expressions asymptotiques
$$
A_z=g^{-1}\left(
\begin{array}{cc} 
   A_{\zeta }                    & -\frac{i}{2}d\;\overline{{\bf jet}_2 B}
\\
\\
\frac{i}{2}d\;{\bf jet}_2 B & \overline{A}_{\zeta } 
\end{array}  
\right)+O(|z|^3)
=\left(
\begin{array}{cc} 
E                   & -\frac{i}{2}d\;\overline{{\bf jet}_2 B}
\\
\\
\frac{i}{2}d\;{\bf jet}_2 B & \overline{E}
\end{array}  
\right)+O(|z|^3),
$$
ce qui nous donne l'expression voulue de la connexion de Chern. Le fait que le repère $(\zeta _k+\bar{\zeta}_k)_k\in{\cal E}(T_X)(U)$ soit $\omega (x)$-orthonormé en $x$ entraîne  qu'on dispose de l'égalité $\eqref{ponctcurvat} $ au point $x$. 
On en déduit donc la formule
$$
C_{m,l}(x)=\Big(\bar{\partial}_{_J }\partial_{_J } h_{l,m}-
\sum_r\,\bar{\partial}_{_J } h_{r,m}\wedge\partial_{_J } h_{l,r} +\partial_{_J }(A''_{\zeta} )_{m,l}-
\bar{\partial}_{_J }\overline{(A''_{\zeta} )_{l,m}}\Big)(x).
$$
pour les coefficients de l'expression locale $\eqref{localcurvattg} $ du tenseur de courbure de Chern du fibré tangent.
On a l'expression
$$
\bar{\partial}_{_J }\partial_{_J } H=-\sum_{j,k}\,\frac{\partial}{\partial \bar{z}_k }\,(\zeta _j\,.H)\,dz_j\wedge d\bar{z}_k   +O(|z|)
$$
avec $\zeta _j\,.H=\sum_p\,(2H^{j,p}\,z_p+H^{j,\bar{p}}\,\bar{z}_p)+O(|z|^2)$. On déduit alors l'expression $$
\bar{\partial}_{_J }\partial_{_J } H=-\sum_{j,k}\, H^{j,\bar{k}}\,dz_j\wedge d\bar{z}_k   +O(|z|).
$$
On a aussi l'expression
$$
\partial_{_J } A''_{\zeta} =\sum_{j,k}\,\frac{\partial}{\partial \bar{z}_j }\,(A''_{\zeta })^k\,dz_j\wedge d\bar{z}_k   +O(|z|).
$$
En rappelant l'expression normale asymptotique $\eqref{conexionA}$ de la forme de connexion $A''_{\zeta}$ par rapport au repère normal 
$(\zeta_k)_k$ on déduit l'expression
$$
\partial_{_J }(A''_{\zeta} )_{m,l}=\frac{1}{4}\sum_{j,k,r}\, (\overline{B}^k_{m,r}-\overline{B}^r_{m,k} )B^j_{r,l} \,dz_j\wedge d\bar{z}_k   +O(|z|).   
$$
En combinant les expressions ainsi obtenues on obtient l'expression $\eqref{valcurvat}$ pour les coefficients $C^{j,k}_{m,l}(x)$ de la courbure au point $x$.\hfill $\Box$
\subsection{Le cas d'une métrique symplectique sur une variété presque complexe}
Dans le cas où la variété presque complexe admet une métrique symplectique, certains des coefficients du lemme précédent se simplifient. On a le lemme suivant.
\begin{lem} 
Soit $(X,J)$ une variété presque complexe admettant une métrique symplectique $\omega\in {\cal E}(\Lambda ^{1,1}_{_J}T_X^*)(X)$. Pour tout point $x$ on peut choisir des coordonnées presque complexes $(z_1,...,z_n)$ d'ordre $N\geq 2$ en $x$ telles que
\begin{eqnarray*}
&\displaystyle{
\omega =\frac{i}{2}\sum_l dz_l\wedge d\bar{z}_l+\frac{i}{2}\sum_{l,m,j,k}\,\Big[H^{j,k}_{l,m}\,z_jz_k+\overline{H}^{j,k}_{m,l}\,\bar{z}_j\bar{z}_k
+\Big(H^{j,\bar{k}}_{l,m}+\frac{i}{4}\sum_r\,B^j_{r,l}\overline{B}^k_{r,m}\Big)z_j\bar{z}_k  \Big]\,dz_l\wedge d\bar{z}_m }&
\\
\\  
&\displaystyle{-\frac{1}{4}\sum_{l,m}\,{\bf jet}_2 B_{l,m}(z)\,dz_l\wedge dz_m
-\frac{1}{4}\sum_{l,m} \,\overline{{\bf jet}_2 B_{l,m}(z)}\,d\bar{z}_l\wedge d\bar{z}_m+O(|z|^3)}.&
\end{eqnarray*}
Quels que soient les coordonnées presque complexes $(z_1,...,z_n)$ d'ordre $N\geq 2$ en $x$  pour lesquelles la métrique $\omega $ s'écrit sous la forme précédente on a l'expression suivante pour le tenseur de courbure de Chern.
\begin{eqnarray*}
&\displaystyle{{\cal C}_{\omega }(T_{X,J})=\sum_{j,k,m,l}\,C^{j,k}_{m,l}(0)\,dz_j\wedge d\bar{z}_k\otimes dz_l
\otimes_{_{J_0} }\frac{\partial}{\partial z_m}    +O(|z|)    }& 
\end{eqnarray*}
avec
\begin{eqnarray*}
C^{j,k}_{m,l}(0)=-H^{j,\bar{k}}_{l,m}+ \frac{1}{4}\sum_r\,\Big[(\overline{B}^k_{m,r} -\overline{B}^r_{m,k})B^j_{r,l} 
+(B^j_{l,r}-B^r_{l,j})\overline{B}^k_{r,m}\Big].
\end{eqnarray*}
\end{lem} 
Contrairement au cas Kählerien, (voir \cite{B-D-I-P}) on ne peut pas éliminer les termes $H^{j,k}_{l,m}\,z_jz_k$ et $\overline{H}^{j,k}_{m,l}\,\bar{z}_j\bar{z}_k$. L'obstruction dérive des termes d'ordre un du jet de la torsion de la structure presque complexe.
\\
\\
$Preuve$
\\
Soient $(z_1,...,z_n)$ des coordonnées presque complexes d'ordre un au point $x$ et $(\zeta _k+\bar{\zeta}_k)_k\in{\cal E}(T_X)(U)$, $\zeta _k=(\partial/\partial z_k)^{1,0}_{_{J} }$ un repère $\omega (x)$-orthonormé. En considérant l'expression du jet d'ordre un du repère $(\zeta ^*_k)_k$ on obtient l'expression locale suivante de la métrique
\begin{eqnarray*}
&\displaystyle{
\omega =\frac{i}{2}\sum_l dz_l\wedge d\bar{z}_l+\frac{i}{2}\sum_{l,m,p}\,\Big(H^p_{l,m}\,z_p+\overline{H}^p_{m,l}\,\bar{z}_p \Big)\,dz_l\wedge d\bar{z}_m }&
\\
\\
&\displaystyle{-\frac{1}{4}\sum_{l,m,p}\, B^p_{l,m}\,z_p\,dz_l\wedge dz_m
-\frac{1}{4}\sum_{l,m} \,\overline{B}^p_{l,m,p}\,\bar{z}_p\,d\bar{z}_l\wedge d\bar{z}_m+O(|z|^2)}.& 
\end{eqnarray*}
Le fait que la métrique $\omega $ soit symplectique implique l'égalité $H^p_{l,m}=H^l_{p,m}$. En effectuant le changement de variables 
$$
Z_m=z_m+\frac{1}{2}\sum_{p,l}\,H^p_{l,m}\,z_pz_l
$$
on obtient, d'après le corollaire $\ref{corholo}$, des coordonnées presque complexes 
$(Z_1,...,Z_n)$ à l'ordre un en $x$ avec les mêmes coefficients $B^{*,*}_{*,*}$ du jet d'ordre un de la structure presque complexe. L'expression de la métrique par rapport aux nouvelles coordonnées est
\begin{eqnarray*}
\omega =\frac{i}{2}\sum_l dZ_l\wedge d\bar{Z}_l-\frac{1}{4}\sum_{l,m,p}\, B^p_{l,m}\,Z_p\,dZ_l\wedge dZ_m
-\frac{1}{4}\sum_{l,m} \,\overline{B}^p_{l,m,p}\,\bar{Z}_p\,d\bar{Z}_l\wedge d\bar{Z}_m+O(|Z|^2).
\end{eqnarray*}
A partir des coordonnées ainsi obtenues on peut construire (d'après la preuve de la proposition $\ref{JN} $) des coordonnées presque complexes d'ordre $N\geq 2$ en $x$ tout en conservant les coefficients $B^{*,*}_{*,*}$ du jet d'ordre un de $J$. En tenant compte de l'expression $\eqref{repdual}$ du jet d'ordre deux du repère $(\zeta ^*_k)_k$ par rapport aux coordonnées en question  on déduit facilement que la métrique $\omega$ s'écrit sous la forme donnée dans l'énoncé du lemme.\hfill$\Box$ 
\subsection{Expression asymptotique normale du flot géodésique d'une connexion de Chern sur le fibré tangent}\label{Exp-asymp-FlotChern} 
On rappelle que par définition $\exp_z(v):=\gamma (1)$, où $\gamma :[0,1]\longrightarrow X$ est la courbe géodésique solution de l'équation différentielle ordinaire $(\gamma ^*D^{\omega} _{_{J} })\dot{\gamma} =0,\,\dot{\gamma}:=d\gamma /dt\in {\cal E}(\gamma ^*T_X)((0,1]) $ avec les conditions initiales $\gamma (0)=z$ et $\dot{\gamma}(0)=v$.
Le résultat suivant est une généralisation dans le cas presque complexe non intégrable d'un calcul fait par Demailly dans \cite{Dem-2}. 
\begin{theo}
Soit $(X,J)$ une variété presque complexe, $\omega\in {\cal E}(\Lambda ^{1,1}_{_J}T_X^*)(X)$ une métrique hermitienne et soient $(z_1,...,z_n)$  des coordonnées presque complexes d'ordre $N\geq 2$ en un point $x$ telles que le repère normal $(\zeta _k+\bar{\zeta}_k)_k\in{\cal E}(T_X)(U_x),\,\zeta _k=(\partial/\partial z_k)^{1,0}_{_{J} }$ soit $\omega (x)$-orthonormé. Le flot géodésique $\exp: {\cal U}\subset T_X\longrightarrow X$  induit par la connexion de Chern du fibré tangent
$$
D^{\omega} _{_{J} }: {\cal E}(T_{X,J}) \longrightarrow {\cal E}(T^*_X\otimes_{_{\R}}T_{X,J})$$
associé à la métrique $\omega$, $($ici ${\cal U}\subset T_X$ désigne un voisinage ouvert de la section nulle$)$, admet l'expression asymptotique suivante au point $(x,0)\in T_{X,x}$ :
\begin{eqnarray*}
&\displaystyle{
\exp_z(v)_k=z_k+v_k-\frac{1}{2}\sum_{l,p,h}\Big[\Big(\hat{S}^{p,h}_{k,l}\,z_h+ S^{p,\bar{h}}_{k,l}\,\bar{z}_h\Big)v_pv_l+\Big(S^{\bar{p},h}_{k,l}\,z_h+S^{\bar{p},\bar{h}}_{k,l}\,\bar{z}_h\Big)\bar{v}_pv_l \Big] }&
\\
\\
&\displaystyle{+\frac{i}{4}\sum_{p,l}\Big[\overline{B}^p_{k,l}+\sum_h\Big(\overline{B}^{p,\bar{h}}_{k,l}\,z_h+2\overline{B}^{p,h}_{k,l}\,\bar{z}_h\Big)\Big]\bar{v}_p\bar{v}_l+O(|v|^2(|z|^2+|v|)) }& 
\end{eqnarray*}  
où 
$
v= \sum_k\,(v_k \frac{\partial}{\partial z_k}+\bar{v}_k\frac{\partial}{\partial \bar{z}_k})\in T_{X,z} 
$, 
$
\hat{S}^{p,h}_{k,l}:=2H^{p,h}_{l,k}-\frac{i}{2}B^{h,\bar{k}}_{l,p}
-\frac{i}{2}\sum_j\overline{H}^j_{k,l}B^h_{j,p}
$
et les autres coefficients $S^{*,*}_{k,l}$ sont donnés dans l'énonce du lemme \ref{Exp-asym-Conn-Chern}.
\end{theo} 
$Preuve$. Par rapport aux coordonnées presque complexes en question nous considérons les écritures $\gamma (t)\equiv(\gamma _1(t),...,\gamma _n(t))$ et
$$
\dot{\gamma}(t)=\sum_k\Big(\dot{\gamma}_k(t)\frac{\partial}{\partial z_k}{\vphantom{z_k} } _{| _{_\gamma (t)} }+\overline{\dot{\gamma}_k(t)} \frac{\partial}{\partial \bar{z} _k}{\vphantom{\bar{z}_k} } _{| _{_\gamma (t)} } \Big).
$$
On pose par définition $\ddot{\gamma}_k:=d^2\gamma _k/dt^2$. On déduit d'après le lemme \ref{Exp-asym-Conn-Chern} que l'équation différentielle ordinaire $(\gamma ^*D^{\omega} _{_{J} })\dot{\gamma} =0$ s'écrit sous la forme
\begin{eqnarray}\label{eq-géodesiq}  
\ddot{\gamma}_k(t)+\sum_l\,\Big[E_{k,l\,|_{\gamma (t)}}(\dot{\gamma}(t))\cdot\dot{\gamma}_l(t)  -
\frac{i}{2}d\;\overline{ {\bf jet}_2B_{k,l\,|_{\gamma (t)}} }(\dot{\gamma}(t))\cdot \overline{\dot{\gamma}_l(t)} \Big]
+O(|\gamma (t)|^2) |\dot{\gamma}(t)|^2=0.\;\;
\end{eqnarray} 
Les conditions initiales $\gamma (0)=z$ et $\dot{\gamma}(0)=v$ donnent l'expression asymptotique $\gamma _k(t)=z_k+tv_k+O(t^2|v|^2)$. On remarque que si $\tau _{_{J}}(x)=0$ alors le terme d'erreur est $O(t^2|z|\, |v|^2)$. En remplaçant l'expression précédente dans l'équation \eqref{eq-géodesiq} et en remarquant qu'on peut toujours supposer la condition $H^p_{l,k}=-H^l_{p,k}$ on obtient l'expression asymptotique suivante pour les dérivées deuxièmes de la courbe $\gamma$ :
\begin{eqnarray*}
&\displaystyle{
\ddot{\gamma}_k(t)=-\sum_{l,p,h}\Big[\Big(\hat{S}^{p,h}_{k,l}\,z_h+ S^{p,\bar{h}}_{k,l}\,\bar{z}_h\Big)v_pv_l+\Big(S^{\bar{p},h}_{k,l}\,z_h+S^{\bar{p},\bar{h}}_{k,l}\,\bar{z}_h\Big)\bar{v}_pv_l \Big] }&
\\
\\
&\displaystyle{+\frac{i}{2}\sum_{p,l}\Big[\overline{B}^p_{k,l}+\sum_h\Big(\overline{B}^{p,\bar{h}}_{k,l}\,z_h+2\overline{B}^{p,h}_{k,l}\,\bar{z}_h\Big)\Big]\bar{v}_p\bar{v}_l+O(|v|^2(|z|^2+|v|))(t)}.& 
\end{eqnarray*}
Si $\tau _{_{J}}(x)=0$ alors le calcul peut être effectue avec plus de précision  car les termes $\overline{B}^p_{k,l}$ sont nuls dans ce cas. Le terme d'erreur serait alors $O(|v|^2(|z|^2+|v|)^2)(t)$.
En intégrant deux fois de suite l'expression précédente on obtient l'expression asymptotique
\begin{eqnarray*}
&\displaystyle{
\gamma_k(t)=z_k+tv_k-\frac{t^2}{2}\sum_{l,p,h}\Big[\Big(\hat{S}^{p,h}_{k,l}\,z_h+ S^{p,\bar{h}}_{k,l}\,\bar{z}_h\Big)v_pv_l+\Big(S^{\bar{p},h}_{k,l}\,z_h+S^{\bar{p},\bar{h}}_{k,l}\,\bar{z}_h\Big)\bar{v}_pv_l \Big] }&
\\
\\
&\displaystyle{+\frac{it^2}{4}\sum_{p,l}\Big[\overline{B}^p_{k,l}+\sum_h\Big(\overline{B}^{p,\bar{h}}_{k,l}\,z_h+2\overline{B}^{p,h}_{k,l}\,\bar{z}_h\Big)\Big]\bar{v}_p\bar{v}_l+O(|v|^2(|z|^2+|v|))(t) }& 
\end{eqnarray*} 
qui permet de conclure la preuve du théorème.\hfill$\Box$

\vspace{1cm}
\noindent
Nefton Pali
\\
Institut Fourier, UMR 5582, Université Joseph Fourier
\\
BP 74, 38402 St-Martin-d'Hères cedex, France
\\
E-mail: \textit{nefton.pali@ujf-grenoble.fr}
\end{document}